\documentclass{article}
\usepackage{amsmath}
\usepackage{amssymb}
\usepackage{amsthm}
\usepackage{cite}
\usepackage[arrow,curve,matrix]{xy}
\begin{document}
\def\e#1\e{\begin{equation}#1\end{equation}}
\def\ea#1\ea{\begin{align}#1\end{align}}
\def\eq#1{{\rm(\ref{#1})}}
\theoremstyle{plain}
\newtheorem{thm}{Theorem}[section]
\newtheorem{prop}[thm]{Proposition}
\newtheorem{cor}[thm]{Corollary}
\theoremstyle{definition}
\newtheorem{dfn}[thm]{Definition}
\newtheorem{ex}[thm]{Example}
\newtheorem{conv}[thm]{Convention}
\def\dim{\mathop{\rm dim}\nolimits}
\def\vdim{\mathop{\rm vdim}\nolimits}
\def\Ker{\mathop{\rm Ker}}
\def\Im{\mathop{\rm Im}}
\def\Kh{{\rm Kh}}
\def\Kb{{\rm Kb}}
\def\Kcb{{\rm Kcb}}
\def\cs{{\rm cs}}
\def\Kh{{\rm Kh}}
\def\Kch{{\rm Kch}}
\def\rsi{{\rm si}}
\def\bo{{\rm bo}}
\def\na{{\rm na}}
\def\Fix{\mathop{\rm Fix}}
\def\Pd{\mathop{\rm Pd}\nolimits}
\def\rank{\mathop{\rm rank}}
\def\Stab{\mathop{\rm Stab}\nolimits}
\def\Aut{\mathop{\rm Aut}}
\def\id{\mathop{\rm id}\nolimits}
\def\bs{\boldsymbol}
\def\ge{\geqslant}
\def\le{\leqslant\nobreak}
\def\N{{\mathbin{\mathbb N}}}
\def\R{{\mathbin{\mathbb R}}}
\def\Z{{\mathbin{\mathbb Z}}}
\def\Q{{\mathbin{\mathbb Q}}}
\def\C{{\mathbin{\mathbb C}}}
\def\oM{{\mathbin{\smash{\,\,\overline{\!\!\mathcal M\!}\,}}}}
\def\al{\alpha}
\def\be{\beta}
\def\ga{\gamma}
\def\de{\delta}
\def\io{\iota}
\def\ep{\epsilon}
\def\la{\lambda}
\def\ka{\kappa}
\def\th{\theta}
\def\si{\sigma}
\def\om{\omega}
\def\De{\Delta}
\def\La{\Lambda}
\def\Om{\Omega}
\def\Ga{\Gamma}
\def\pd{\partial}
\def\ts{\textstyle}
\def\op{\oplus}
\def\ot{\otimes}
\def\iy{\infty}
\def\ra{\rightarrow}
\def\ab{\allowbreak}
\def\longra{\longrightarrow}
\def\t{\times}
\def\ci{\circ}
\def\ti{\tilde}
\def\es{\emptyset}
\def\d{{\rm d}}
\def\ha{{\ts\frac{1}{2}}}
\def\md#1{\vert #1 \vert}
\def\bmd#1{\big\vert #1 \big\vert}
\title{Kuranishi homology and Kuranishi cohomology:\\
a User's Guide}
\author{Dominic Joyce}
\date{}
\maketitle
\begin{abstract}
A {\it Kuranishi space\/} is a topological space equipped with a
{\it Kuranishi structure}, defined by Fukaya and Ono. Kuranishi
structures occur naturally on many moduli spaces in differential
geometry, and in particular, in moduli spaces of stable
$J$-holomorphic curves in Symplectic Geometry.

This paper is a summary of the author's book \cite{Joyc1}. Let $Y$
be an orbifold, and $R$ a $\Q$-algebra. We shall define a new
homology theory of $Y$, {\it Kuranishi homology} $KH_*(Y;R)$, using
a chain complex $KC_*(Y;R)$ spanned by isomorphism classes $[X,\bs
f,\bs G]$, where $X$ is a compact, oriented Kuranishi space with
corners, $\bs f:X\ra Y$ is strongly smooth, and $\bs G$ is some
extra {\it gauge-fixing data\/} for $(X,\bs f)$. The purpose of $\bs
G$ is to ensure the automorphism groups $\Aut(X,\bs f,\bs G)$ are
finite, which is necessary to get a well-behaved homology theory.
Our main result is that $KH_*(Y;R)$ is isomorphic to singular
homology.

We define Poincar\'e dual {\it Kuranishi (co)homology} $KH^*(Y;R)$,
which is isomorphic to compactly-supported cohomology, using a
cochain complex $KC^*(Y;R)$ spanned by $[X,\bs f,\bs C]$, where $X$
is a compact Kuranishi space with corners, $\bs f:X\ra Y$ is a
cooriented strong submersion, and $\bs C$ is {\it co-gauge-fixing
data}. We also define simpler theories of {\it Kuranishi bordism\/}
and {\it Kuranishi cobordism} $KB_*,KB^*(Y;R)$, for $R$ a
commutative ring.

These theories are powerful new tools in Symplectic Geometry. Moduli
spaces of $J$-holomorphic curves define (co)chains directly in
Kuranishi (co)homology or Kuranishi (co)bordism. This hugely
simplifies the formation of virtual cycles, as there is no longer
any need to perturb moduli spaces. The theory has applications to
Lagrangian Floer cohomology, String Topology, and the
Gopakumar--Vafa Integrality Conjecture.
\end{abstract}

\tableofcontents

\section{Introduction}
\label{ug1}

A {\it Kuranishi space\/} is a topological space with a {\it
Kuranishi structure}, defined by Fukaya and Ono \cite{FuOn,FOOO}.
Let $Y$ be an orbifold and $R$ a $\Q$-algebra. In the book
\cite{Joyc1} the author develops {\it Kuranishi homology\/}
$KH_*(Y;R)$ and {\it Kuranishi cohomology\/} $KH^*(Y;R)$. The
(co)chains in these theories are of the form $[X,\bs f,\bs G]$ where
$X$ is a compact Kuranishi space, $\bs f:X\ra Y$ is a strongly
smooth map, and $\bs G$ is some extra {\it gauge-fixing data}. We
prove $KH_*(Y;R)$ is isomorphic to singular homology
$H_*^\rsi(Y;R)$, and $KH^*(Y;R)$ is isomorphic to
compactly-supported cohomology $H^*_\cs(Y;R)$. We also define {\it
Kuranishi bordism\/} $KB_*(Y;R)$ and {\it Kuranishi cobordism}
$KB^*(Y;R)$, for $R$ a commutative ring.

This paper is a brief introduction to selected parts of
\cite{Joyc1}. The length of \cite{Joyc1} (presently 290 pages) is
likely to deter people from reading it, but the main ideas can be
summarized much more briefly, and that is what this paper tries to
do. A User's Guide, say for a car or a computer, should give you a
broad overview of how the machine actually works, and instructions
on how to use it in practice, but it should probably not tell you
exactly where the carburettor is, or how the motherboard is wired.
This paper is written in the same spirit. It should provide you with
sufficient background to understand the sequels
\cite{AkJo,Joyc2,Joyc3} (once I get round to writing them), and also
to decide whether Kuranishi (co)homology is a good tool to use in
problems you are interested in.

Here are the main areas of \cite{Joyc1} that we will {\it not\/}
cover. In \cite{Joyc1} we also define {\it effective Kuranishi
(co)homology\/} $KH_*^{\rm ef},KH^*_{\rm ec}(Y;R)$, a variant on
Kuranishi homology that has the advantage that it works for all
commutative rings $R$, including $R=\Z$, and is isomorphic to
$H_*^\rsi,H^*_\cs(Y;R)$, but has the disadvantage of restrictions on
the Kuranishi spaces $X$ allowed in chains $[X,\bs f,\bs G]$, which
limits its applications in Symplectic Geometry. Similarly,
\cite[Ch.~5]{Joyc1} actually defines five different kinds of
Kuranishi (co)bordism, but below we consider only the simplest. The
applications to Symplectic Geometry in \cite[Ch.~6]{Joyc1} are
omitted.

Kuranishi (co)homology is intended primarily as a tool for use in
areas of Symplectic Geometry involving $J$-holomorphic curves, and
will be applied in the sequels \cite{AkJo,Joyc2,Joyc3,Joyc4}.
Kuranishi structures occur naturally on many moduli spaces in
Differential Geometry. For example, if $(M,\om)$ is a compact
symplectic manifold with almost complex structure $J$ then the
moduli space $\oM_{g,m}(M,J,\be)$ of stable $J$-holomorphic curves
of genus $g$ with $m$ marked points in class $\be$ in $H_2(M;\Z)$ is
a compact Kuranishi space with a strongly smooth map $\prod_i{\bf
ev}_i:\oM_{g,m}(M,J,\be)\ra M^m$. By choosing some gauge-fixing data
$\bs G$ we define a cycle $[\oM_{g,m}(M,J,\be),\prod_i{\bf ev}_i,\bs
G]$ in $KC_*(M^m;\Q)$ whose homology class $\bigl[[\oM_{g,m}(M,\ab
J,\ab\be),\prod_i{\bf ev}_i,\bs G]\bigr]$ in $KH_*(M^m;\Q)\cong
H^\rsi_*(M^m;\Q)$ is a {\it Gromov--Witten invariant\/} of
$(M,\om)$, and is independent of the choice of almost complex
structure~$J$.

In the conventional definitions of symplectic Gromov--Witten
invariants \cite{FuOn,LiTi,Ruan,Sieb}, one must define a {\it
virtual cycle\/} for $\oM_{g,m}(M,J,\be)$. This is a complicated
process, involving many arbitrary choices: first one must perturb
the moduli space, morally over $\Q$ rather than $\Z$, so that it
becomes something like a manifold. Then one must triangulate the
perturbed moduli space by simplices to define a cycle in the
singular chains $C_*^\rsi(M^m;\Q)$. The Gromov--Witten invariant is
the homology class of this virtual cycle. By using Kuranishi
(co)homology as a substitute for singular homology, this process of
defining virtual cycles becomes much simpler and less arbitrary.
{\it The moduli space is its own virtual cycle}, and we eliminate
the need to perturb moduli spaces and triangulate by simplices.

The real benefits of the Kuranishi (co)homology approach come not in
closed Gromov--Witten theory, where the moduli spaces are Kuranishi
spaces without boundary, but in areas such as open Gromov--Witten
theory, Lagrangian Floer cohomology \cite{FOOO}, Contact Homology
\cite{EES}, and Symplectic Field Theory \cite{EGH}, where the moduli
spaces are Kuranishi spaces {\it with boundary and corners}, and
their boundaries are identified with fibre products of other moduli
spaces.

In the conventional approach, one must choose virtual chains for
each moduli space, which must be compatible at the boundary with
intersection products of choices of virtual chains for other moduli
spaces. This business of boundary compatibility of virtual chains is
horribly complicated and messy, and a large part of the 1385 pages
of Fukaya, Oh, Ohta and Ono's work on Lagrangian Floer cohomology
\cite{FOOO} is devoted to dealing with it. Using Kuranishi
cohomology, because we do not perturb moduli spaces, choosing
virtual chains with boundary compatibility is very easy, and
Lagrangian Floer cohomology can be reformulated in a much more
economical way, as we will show in \cite[\S 6.6]{Joyc1}
and~\cite{AkJo}.

An important feature of these theories is that {\it Kuranishi
homology and cohomology are very well behaved at the (co)chain
level}, much better than singular homology, say. For example,
Kuranishi cochains $KC^*(Y;R)$ have a supercommutative, associative
cup product $\cup$, cap products also work well at the (co)chain
level, and there is a well-behaved functor from singular chains
$C_*^\rsi(Y;R)$ to Kuranishi chains $KC_*(Y;R)$. Because of this,
the theories may also have applications in other areas which may not
be directly related to Kuranishi spaces, but which need a
(co)homology theory of manifolds or orbifolds with good
(co)chain-level behaviour. In \cite{Joyc3} we will apply Kuranishi
(co)chains to reformulate the String Topology of Chas and Sullivan
\cite{ChSu}, which involves chains on infinite-dimensional loop
spaces. Another possible area is Costello's approach to Topological
Conformal Field Theories \cite{Cost}, which involves a choice of
chain complex for homology, applied to moduli spaces of Riemann
surfaces.

As well as Kuranishi homology and cohomology, in \cite[Ch.~5]{Joyc1}
we also define {\it Kuranishi bordism\/} $KB_*(Y;R)$ and {\it
Kuranishi cobordism\/} $KB^*(Y;R)$. These are simpler than Kuranishi
(co)homology, being spanned by $[X,\bs f]$ for $X$ a compact
Kuranishi space {\it without boundary} and $\bs f:X\ra Y$ strongly
smooth, and do not involve gauge-fixing data. In contrast to
Kuranishi (co)homology which is isomorphic to conventional homology
and compactly-supported cohomology, these are new topological
invariants, and we show that they are very large --- for instance,
if $Y\ne\es$ and $R\ot_\Z\Q\ne 0$ then $KB_{2k}(Y;R)$ is infinitely
generated over $R$ for all~$k\in\Z$.

In \cite[\S 6.2]{Joyc1} we define new Gromov--Witten type invariants
$[\oM_{g,m}(M,\ab J,\ab\be),\prod_i{\bf ev}_i]$ in Kuranishi bordism
$KB_*(M^m;\Z)$. Since these are defined in groups over $\Z$, not
$\Q$, the author expects that Kuranishi (co)bordism will be useful
in studying {\it integrality properties} of Gromov--Witten
invariants. In \cite[\S 6.3]{Joyc1} we outline an approach to
proving the Gopakumar--Vafa Integrality Conjecture for
Gromov--Witten invariants of Calabi--Yau 3-folds, which the author
hopes to take further in~\cite{Joyc4}.
\medskip

\noindent{\it Acknowledgements.} I am grateful to Mohammed Abouzaid,
Manabu Akaho, Kenji Fukaya, Ezra Getzler, Shinichiroh Matsuo,
Yong-Geun Oh, Hiroshi Ohta, Kauru Ono, Paul Seidel, Ivan Smith and
Dennis Sullivan for useful conversations. This research was
partially supported by EPSRC grant EP/D07763X/1.

\section{Kuranishi spaces}
\label{ug2}

{\it Kuranishi spaces\/} were introduced by Fukaya and Ono
\cite{FOOO,FuOn}, and are important in Symplectic Geometry because
moduli spaces of stable $J$-holomorphic curves in symplectic
manifolds are Kuranishi spaces. We use the definitions and notation
of \cite[\S 2]{Joyc1}, which have some modifications from those of
Fukaya--Ono.

\subsection{Manifolds and orbifolds with corners and g-corners}
\label{ug21}

In \cite{Joyc1} we work with four classes of manifolds, in
increasing order of generality: {\it manifolds without boundary},
{\it manifolds with boundary}, {\it manifolds with corners}, and
{\it manifolds with generalized corners\/} or {\it g-corners}. The
first three classes are fairly standard, although the author has not
found a reference for foundational material on manifolds with
corners. Manifolds with g-corners are new, as far as the author
knows. The precise definitions of these classes of manifolds are
given in \cite[\S 2.1]{Joyc1}. Here are the basic ideas:
\begin{itemize}
\setlength{\itemsep}{0pt}
\setlength{\parsep}{0pt}
\item An $n$-{\it dimensional manifold without boundary\/} is
locally modelled on open sets in $\R^n$.
\item An $n$-{\it dimensional manifold with boundary\/} is
locally modelled on open sets in $\R^n$ or $[0,\iy)\t\R^{n-1}$.
\item An $n$-{\it dimensional manifold with corners\/} is
locally modelled on open sets in $[0,\iy)^k\t\R^{n-k}$ for
$k=0,\ldots,n$.
\item A {\it polyhedral cone\/} $C$ in $\R^n$ is a subset of the
form
\begin{equation*}
C=\bigl\{(x_1,\ldots,x_n)\in\R^n: a_1^ix_1+\cdots+a_n^ix_n\ge 0,\;\>
i=1,\ldots,k\bigr\},
\end{equation*}
where $a^i_j\in\R$ for $i=1,\ldots,k$ and~$j=1,\ldots,n$.

An $n$-{\it dimensional manifold with g-corners\/} is roughly
speaking locally modelled on open sets in polyhedral cones $C$ in
$\R^n$ with nonempty interior $C^\ci$. Since $[0,\iy)^k\t\R^{n-k}$
is a polyhedral cone $C$ with $C^\ci\ne\es$, manifolds with corners
are examples of manifolds with g-corners. In fact the subsets in
$\R^n$ used as local models for manifolds with g-corners are more
general than polyhedral cones, but this extra generality is only
needed for technical reasons in the proof of Theorem \ref{ug3thm2}.
\end{itemize}

Here are some examples. The line $\R$ is a 1-manifold without
boundary; the interval $[0,1]$ is a 1-manifold with boundary; the
square $[0,1]^2$ is a 2-manifold with corners; and the octahedron
\begin{equation*}
O=\bigl\{(x_1,x_2,x_3)\in\R^3:\md{x_1}+\md{x_2}+\md{x_3}\le 1\bigr\}
\end{equation*}
in $\R^3$ is a 3-manifold with g-corners. It is not a manifold with
corners, since four 2-dimensional faces of $O$ meet at the vertex
$(1,0,0)$, but three 2-dimensional faces of $[0,\iy)^3$ meet at the
vertex $(0,0,0)$, so $O$ near $(1,0,0)$ is not locally modelled on
$[0,\iy)^3$ near~$(0,0,0)$.

Manifolds $X$ with boundary, corners, or g-corners have a
well-behaved notion of {\it boundary\/} $\pd X$. To motivate the
definition, consider $[0,\iy)^2$ in $\R^2$. If we took
$\pd\bigl([0,\iy)^2\bigr)$ to be the subset
$\bigl([0,\iy)\t\{0\}\bigr)\cup\bigl(\{0\}\t[0,\iy)\bigr)$ of
$[0,\iy)^2$, then $\pd\bigl([0,\iy)^2\bigr)$ would not be a manifold
with corners near $(0,0)$. Instead, we take
$\pd\bigl([0,\iy)^2\bigr)$ to be the {\it disjoint union\/} of the
two boundary strata $[0,\iy)\t\{0\}$ and $\{0\}\t[0,\iy)$. This is a
manifold with boundary, but now $\pd\bigl([0,\iy)^2\bigr)$ is {\it
not a subset of\/} $[0,\iy)^2$, since two points in
$\pd\bigl([0,\iy)^2\bigr)$ correspond to $(0,0)$ in~$[0,\iy)^2$.

We define the {\it boundary\/} $\pd X$ of an $n$-manifold $X$ with
(g-)corners to be the set of pairs $(p,B)$, where $p\in X$ and $B$
is a local choice of connected $(n\!-\!1)$-dimensional boundary
stratum of $X$ containing $p$. Thus, if $p$ lies in a codimension
$k$ corner of $X$ locally modelled on $[0,\iy)^k\t\R^{n-k}$ then $p$
is represented by $k$ distinct points $(p,B_i)$ in $\pd X$ for
$i=1,\ldots,k$. Then $\pd X$ is an $(n\!-\!1)$-manifold with
(g-)corners. Note that $\pd X$ is not a subset of $X$, but has a
natural immersion $\io:\pd X\ra X$ mapping $(p,B)\mapsto p$. Often
we suppress $\io$, and talk of restricting data on $X$ to $\pd X$,
when really we mean the pullback by~$\io$.

If $X$ is a $n$-manifold with (g-)corners then $\pd^2X$ is an
$(n\!-\!2)$-manifold with (g-)corners. Points of $\pd^2X$ may be
written $(p,B_1,B_2)$, where $p\in X$ and $B_1,B_2$ are distinct
local boundary components of $X$ containing $p$. There is a natural,
free involution $\si:\pd^2X\ra\pd^2X$ mapping
$\si:(p,B_1,B_2)\mapsto(p,B_2,B_1)$, which is orientation-reversing
if $X$ is oriented. This involution is important in questions to do
with extending data defined on $\pd X$ to $X$. For example, if
$f:\pd X\ra\R$ is a smooth function, then a necessary condition for
there to exist smooth $g:X\ra\R$ with $g\vert_{\pd X}\equiv f$ is
that $f\vert_{\pd^2X}$ is $\si$-invariant, and if $X$ has corners
(not g-corners) then this condition is also sufficient.

{\it Orbifolds\/} are a generalization of manifold, which allow
quotients by finite groups. Again, we define orbifolds {\it without
boundary}, or {\it with boundary}, or {\it with corners}, or {\it
with g-corners}, where orbifolds without boundary are locally
modelled on quotients $\R^n/\Ga$ for $\Ga$ a finite group acting
linearly on $\R^n$, and similarly for the other classes. (We do not
require $\Ga$ to act {\it effectively}, so we cannot regard $\Ga$ as
a subgroup of GL$(n,\R)$.) Orbifolds (then called $V$-manifolds)
were introduced by Satake \cite{Sata}, and a book on orbifolds is
Adem et al.\ \cite{ALR}. Note however that the right definition of
smooth maps of orbifolds is not that given by Satake, but the more
complex notion of {\it morphisms of orbifolds\/} in \cite[\S
2.4]{ALR}. When we do not specify otherwise, by a manifold or
orbifold, we always mean a manifold or orbifold with g-corners, the
most general class.

Let $X,Y$ be manifolds of dimensions $m,n$. {\it Smooth maps}
$f:X\ra Y$ are continuous maps which are locally modelled on smooth
maps from $\R^m\ra\R^n$. A smooth map $f$ induces a morphism of
vector bundles $\d f:TX\ra f^*(TY)$ on $X$, where $TX,TY$ are the
tangent bundles of $X$ and $Y$. For manifolds $X,Y$ without
boundary, we call a smooth map $f:X\ra Y$ a {\it submersion} if $\d
f:TX\ra f^*(TY)$ is a surjective morphism of vector bundles. If
$X,Y$ have boundary or (g-)corners, the definition of submersions
$f:X\ra Y$ in \cite[\S 2.1]{Joyc1} is more complicated, involving
conditions over $\pd^kX$ and $\pd^lY$ for all $k,l\ge 0$. When $\pd
Y=\es$, $f$ is a submersion if $\d(f\vert_{\pd^kX}):T(\pd^kX)\ra
f\vert_{\pd^kX}^*(TY)$ is surjective for all~$k\ge 0$.

Let $X,X',Y$ be manifolds (in any of the four classes above) and
$f:X\ra Y$, $f':X'\ra Y$ be smooth maps, at least one of which is a
submersion. Then the {\it fibre product\/} $X\t_{f,Y,f'}X'$ or
$X\t_YX'$ is
\e
X\t_{f,Y,f'}X'=\bigl\{(p,p')\in X\t X':f(p)=f'(p')\bigr\}.
\label{ug2eq1}
\e
It turns out \cite[Prop.~2.6]{Joyc1} that $X\t_YX'$ is a submanifold
of $X\t X'$, and so is a manifold (in the same class as $X,X',Y$).
When $X,X',Y$ have boundary or (g-)corners, the complicated
definition of submersion is necessary to make $X\t_YX'$ a
submanifold over~$\pd^kX,\pd^{k'}X',\pd^lY$.

Fibre products can be defined for orbifolds \cite[\S 2.2]{Joyc1},
but there are some subtleties to do with stabilizer groups. To
explain this, first consider the case in which $U,U',V$ are
manifolds, and $\Ga,\Ga',\De$ are finite groups acting on $U,U',V$
by diffeomorphisms so that $U/\Ga$, $U'/\Ga'$, $V/\De$ are
orbifolds, and $\rho:\Ga\ra\De$, $\rho':\Ga'\ra\De$ are group
homomorphisms, and $f:U\ra V$, $f':U'\ra V$ are smooth $\rho$- and
$\rho'$-equivariant maps, at least one of which is a submersion.
Then $f,f'$ induce smooth maps of orbifolds $f_*:U/\Ga\ra V/\De$,
$f'_*:U'/\Ga'\ra V/\De$, at least one of which is a submersion.

It turns out that the right answer for the orbifold fibre product is
\e
(U/\Ga)\t_{f_*,V/\De,f'_*}(U'/\Ga')=\bigl((U\t U')\t_{f\t f',V\t
V,\pi}(V\t\De)\bigr)/(\Ga\t\Ga').
\label{ug2eq2}
\e
Here $\pi:V\t\De\ra V\t V$ is given by $\pi:(v,\de)\mapsto
(v,\de\cdot v)$, and $(U\t U')\t_{f\t f',V\t V,\pi}(V\t\De)$ is the
fibre product of smooth manifolds, and $\Ga\t\Ga'$ acts on the
manifold $(U\t U')\t_{V\t V}(V\t\De)$ by diffeomorphism
$(\ga,\ga'):\bigl((u,u'),(v,\de)\bigr) \mapsto\bigl((\ga\cdot
u,\ga'\cdot u'),(\rho(\ga)\cdot
v,\rho'(\ga')\de\rho'(\ga')^{-1})\bigr)$, so that the quotient is an
orbifold. Now \eq{ug2eq2} coincides with \eq{ug2eq1} for $X=U/\Ga$,
$X'=U'/\Ga'$, $Y=V/\De$ only if one of $\rho:\Ga\ra\De$,
$\rho',\Ga'\ra\De$ are surjective; otherwise the projection from
\eq{ug2eq2} to \eq{ug2eq1} is a finite surjective map, but not
necessarily injective.

This motivates the definition of fibre products of orbifolds. Let
$X,X',Y$ be orbifolds, and $f:X\ra Y$, $f':X'\ra Y$ be smooth maps,
at least one of which is a submersion. Then for $p\in X$ and $p'\in
X'$ with $f(p)=q=f(p')$ in $Y$ we have morphisms of stabilizer
groups $f_*:\Stab(p)\ra\Stab(q)$, $f'_*:\Stab(p')\ra\Stab(q)$. Thus
we can form the double coset space
\begin{align*}
&f_*(\Stab(p))\backslash\Stab(q)/f'_*(\Stab(p'))\\
&=\bigl\{ \{f_*(\ga)\de f_*(\ga'):\ga\in\Stab(p),\;\>
\ga'\in\Stab(p')\}:\de\in\Stab(q)\bigr\}.
\end{align*}
As a set, we define
\e
\begin{split}
X\t_{f,Y,f'}X'=\bigl\{(p,p',\De):\,&\text{$p\in X$, $p'\in X'$,
$f(p)=f'(p')$,}\\
&\De\in f_*(\Stab(p))\backslash\Stab(f(p))/f'_*(\Stab(p'))\bigr\}.
\end{split}
\label{ug2eq3}
\e
We give this an orbifold structure in a natural way, such that if
$(U,\Ga,\phi)$, $(U',\Ga',\phi'),(V',\De',\psi')$ are orbifold
charts on $X,X',Y$ with $f(\Im\phi),f'(\Im\phi')\subseteq\Im\psi$
then we use \eq{ug2eq2} to define an orbifold chart on~$X\t_YX'$.

\subsection{Kuranishi structures on topological spaces}
\label{ug22}

Let $X$ be a paracompact Hausdorff topological space throughout.

\begin{dfn} A {\it Kuranishi neighbourhood\/} $(V_p,E_p,s_p,
\psi_p)$ of $p\in X$ satisfies:
\begin{itemize}
\setlength{\itemsep}{0pt}
\setlength{\parsep}{0pt}
\item[(i)] $V_p$ is an orbifold, which may
or may not have boundary or (g-)corners;
\item[(ii)] $E_p\ra V_p$ is an orbifold vector bundle over $V_p$,
the {\it obstruction bundle};
\item[(iii)] $s_p:V_p\ra E_p$ is a smooth section, the {\it
Kuranishi map}; and
\item[(iv)] $\psi_p$ is a homeomorphism from $s_p^{-1}(0)$ to
an open neighbourhood of $p$ in $X$, where $s_p^{-1}(0)$ is the
subset of $V_p$ where the section $s_p$ is zero.
\end{itemize}
\label{ug2def1}
\end{dfn}

\begin{dfn} Let $(V_p,E_p,s_p,\psi_p),(\ti V_p,\ti E_p,\ti s_p,
\ti\psi_p)$ be two Kuranishi neighbourhoods of $p\in X$. We call
$(\al,\hat\al):(V_p,\ldots,\psi_p)\ra(\ti V_p,\ldots,\ti\psi_p)$ an
{\it isomorphism\/} if $\al:V_p\ra\ti V_p$ is a diffeomorphism and
$\hat\al:E_p\ra\al^*(\ti E_p)$ an isomorphism of orbibundles, such
that $\ti s_p\ci\al\equiv\hat\al\ci s_p$
and~$\ti\psi_p\ci\al\equiv\psi_p$.

We call $(V_p,\ldots,\psi_p),(\ti V_p,\ldots,\ti\psi_p)$ {\it
equivalent\/} if there exist open neighbourhoods $U_p\!\subseteq\!
V_p$, $\ti U_p\!\subseteq\!\ti V_p$ of $\psi_p^{-1}(p),
\ti\psi_p^{-1}(p)$ such that
$(U_p,E_p\vert_{U_p},s_p\vert_{U_p},\psi_p\vert_{U_p})$ and $(\ti
U_p,\ti E_p\vert_{\ti U_p},\ti s_p\vert_{\ti U_p},\ti\psi_p
\vert_{\ti U_p})$ are isomorphic.
\label{ug2def2}
\end{dfn}

\begin{dfn} Let $(V_p,E_p,s_p,\psi_p)$ and $(V_q,E_q,s_q,
\psi_q)$ be Kuranishi neighbourhoods of $p\in X$ and
$q\in\psi_p(s_p^{-1}(0))$ respectively. We call a pair
$(\phi_{pq},\hat\phi_{pq})$ a {\it coordinate change\/} from
$(V_q,\ldots,\psi_q)$ to $(V_p,\ldots,\psi_p)$ if:
\begin{itemize}
\setlength{\itemsep}{0pt}
\setlength{\parsep}{0pt}
\item[(a)] $\phi_{pq}:V_q\ra V_p$ is a smooth embedding of orbifolds;
\item[(b)] $\hat\phi_{pq}:E_q\ra\phi_{pq}^*(E_p)$ is an embedding of
orbibundles over~$V_q$;
\item[(c)] $\hat\phi_{pq}\ci s_q\equiv s_p\ci\phi_{pq}$;
\item[(d)] $\psi_q\equiv \psi_p\ci\phi_{pq}$; and
\item[(e)] Choose an open neighbourhood $W_{pq}$ of
$\phi_{pq}(V_q)$ in $V_p$, and an orbifold vector subbundle $F_{pq}$
of $E_p\vert_{W_{pq}}$ with $\phi_{pq}^*(F_{pq})=\hat\phi_{pq}
(E_q)$, as orbifold vector subbundles of $\phi_{pq}^*(E_p)$ over
$V_q$. Write $\hat s_p:W_{pq}\ra E_p/F_{pq}$ for the projection of
$s_p\vert_{W_{pq}}$ to the quotient bundle $E_p/F_{pq}$. Now
$s_p\vert_{\phi_{pq}(V_q)}$ lies in $F_{pq}$ by (c), so $\hat
s_p\vert_{\phi_{pq}(V_q)}\equiv 0$. Thus there is a well-defined
derivative
\begin{equation*}
\d\hat s_p:N_{\phi_{pq}(V_q)}V_p\ra
(E_p/F_{pq})\vert_{\phi_{pq}(V_q)},
\end{equation*}
where $N_{\phi_{pq}(V_q)}V_p$ is the normal orbifold vector bundle
of $\phi_{pq}(V_q)$ in $V_p$. Pulling back to $V_q$ using
$\phi_{pq}$, and noting that $\phi_{pq}^*(F_{pq})=\hat\phi_{pq}
(E_q)$, gives a morphism of orbifold vector bundles over~$V_q$:
\e
\d\hat s_p:\frac{\phi_{pq}^*(TV_p)}{(\d\phi_{pq})(TV_q)}\longra
\frac{\phi_{pq}^*(E_p)}{\hat\phi_{pq}(E_q)}\,.
\label{ug2eq4}
\e
We require that \eq{ug2eq4} should be an {\it isomorphism\/} over
$s_q^{-1}(0)$.
\end{itemize}
\label{ug2def3}
\end{dfn}

Here Definition \ref{ug2def3}(e) replaces the notion in
\cite[Def.~5.6]{FuOn}, \cite[Def.~A1.14]{FOOO} that a Kuranishi
structure {\it has a tangent bundle}.

\begin{dfn} A {\it germ of Kuranishi neighbourhoods at\/} $p\in X$
is an equivalence class of Kuranishi neighbourhoods
$(V_p,E_p,s_p,\psi_p)$ of $p$, using the notion of equivalence in
Definition \ref{ug2def2}. Suppose $(V_p,E_p,s_p,\psi_p)$ lies in
such a germ. Then for any open neighbourhood $U_p$ of
$\psi_p^{-1}(p)$ in $V_p$,
$(U_p,E_p\vert_{U_p},s_p\vert_{U_p},\psi_p\vert_{U_p})$ also lies in
the germ. As a shorthand, we say that some condition on the germ
{\it holds for sufficiently small\/} $(V_p,\ldots,\psi_p)$ if
whenever $(V_p,\ldots,\psi_p)$ lies in the germ, the condition holds
for $(U_p,\ldots,\psi_p\vert_{U_p})$ for all sufficiently small
$U_p$ as above.

A {\it Kuranishi structure\/} $\ka$ on $X$ assigns a germ of
Kuranishi neighbourhoods for each $p\in X$ and a {\it germ of
coordinate changes\/} between them in the following sense: for each
$p\in X$, for all sufficiently small $(V_p,\ldots,\psi_p)$ in the
germ at $p$, for all $q\in\Im\psi_p$, and for all sufficiently small
$(V_q,\ldots,\psi_q)$ in the germ at $q$, we are given a coordinate
change $(\phi_{pq},\hat\phi_{pq})$ from $(V_q,\ldots,\psi_q)$ to
$(V_p,\ldots,\psi_p)$. These coordinate changes should be compatible
with equivalence in the germs at $p,q$ in the obvious way, and
satisfy:
\begin{itemize}
\setlength{\itemsep}{0pt}
\setlength{\parsep}{0pt}
\item[(i)] $\dim V_p-\rank E_p$ is independent of $p$\/; and
\item[(ii)] if $q\in\Im\psi_p$ and $r\in\Im\psi_q$ then
$\phi_{pq}\ci\phi_{qr}=\phi_{pr}$ and~$\hat\phi_{pq}\ci
\hat\phi_{qr}=\hat\phi_{pr}$.
\end{itemize}
We call $\vdim X=\dim V_p-\rank E_p$ the {\it virtual dimension\/}
of the Kuranishi structure. A {\it Kuranishi space\/} $(X,\ka)$ is a
topological space $X$ with a Kuranishi structure $\ka$. Usually we
refer to $X$ as the Kuranishi space, suppressing~$\ka$.
\label{ug2def4}
\end{dfn}

Loosely speaking, the above definitions mean that a Kuranishi space
is locally modelled on the zeroes of a smooth section of an orbifold
vector bundle over an orbifold. Moduli spaces of $J$-holomorphic
curves in Symplectic Geometry can be given Kuranishi structures in a
natural way, as in~\cite{FuOn,FOOO}.

\subsection{Strongly smooth maps and strong submersions}
\label{ug23}

In \cite[Def.~2.24]{Joyc1} we define {\it strongly smooth maps\/}
$\bs f:X\ra Y$, for $Y$ an orbifold.

\begin{dfn} Let $X$ be a Kuranishi space, and $Y$ a smooth
orbifold. A {\it strongly smooth map\/} $\bs f:X\ra Y$ consists of,
for all $p\in X$ and all sufficiently small $(V_p,E_p,s_p,\psi_p)$
in the germ of Kuranishi neighbourhoods at $p$, a choice of smooth
map $f_p:V_p\ra Y$, such that for all $q\in\Im\psi_p$ and
sufficiently small $(V_q,\ldots,\psi_q)$ in the germ at $q$ with
coordinate change $(\phi_{pq},\hat\phi_{pq})$ from
$(V_q,\ldots,\psi_q)$ to $(V_p,\ldots,\psi_p)$ in the germ of
coordinate changes, we have $f_p\ci\phi_{pq}=f_q$. Then $\bs f$
induces a continuous map $f:X\ra Y$ in the obvious way.

We call $\bs f$ a {\it strong submersion\/} if all the $f_p$ are
submersions, that is, the maps $\d f_p:TV_p\ra f_p^*(TY)$ are
surjective, and also when $V_p$ has boundary or corners,
$f_p\vert_{\pd V_p}:\pd V_p\ra Y$ is a submersion, and the
restriction of $f_p$ to each codimension $k$ corner is a submersion
for all~$k$.
\label{ug2def5}
\end{dfn}

There is also \cite[Def.~2.25]{Joyc1} a definition of strongly
smooth maps $\bs f:X\ra Y$ for $X,Y$ Kuranishi spaces, which we will
not give. A {\it strong diffeomorphism\/} $\bs f:X\ra Y$ is a
strongly smooth map with a strongly smooth inverse. It is the
natural notion of isomorphism of Kuranishi spaces.

\subsection{Boundaries of Kuranishi spaces}
\label{ug24}

We define the {\it boundary\/} $\pd X$ of a Kuranishi space $X$,
which is itself a Kuranishi space of dimension~$\vdim X-1$.

\begin{dfn} Let $X$ be a Kuranishi space. We shall define a
Kuranishi space $\pd X$ called the {\it boundary\/} of $X$. The
points of $\pd X$ are equivalence classes
$[p,(V_p,\ldots,\psi_p),B]$ of triples $(p,(V_p,\ldots,\psi_p),B)$,
where $p\in X$, $(V_p,\ldots,\psi_p)$ lies in the germ of Kuranishi
neighbourhoods at $p$, and $B$ is a local boundary component of
$V_p$ at $\psi_p^{-1}(p)$. Two triples
$(p,(V_p,\ldots,\psi_p),B),(q,(\ti V_q,\ldots,\ti\psi_q),C)$ are
{\it equivalent\/} if $p=q$, and the Kuranishi neighbourhoods
$(V_p,\ldots,\psi_p),\smash{(\ti V_q,\ldots,\ti\psi_q)}$ are
equivalent so that we are given an isomorphism
$(\al,\hat\al):(U_p,\ldots, \psi_p\vert_{U_p})\ra(\ti U_q,\ldots,
\ti\psi_q\vert_{\ti U_q})$ for open $\psi_p^{-1}(p)\in U_p\subseteq
V_p$ and $\ti\psi_q^{-1}(q)\in\ti U_q\subseteq\ti V_q$, and
$\al_*(B)=C$ near $\ti\psi_q^{-1}(q)$.

We can define a unique natural topology and Kuranishi structure on
$\pd X$, such that $(\pd V_p,E_p\vert_{\pd V_p},s_p\vert_{\pd
V_p},\psi'_p)$ is a Kuranishi neighbourhood on $\pd X$ for each
Kuranishi neighbourhood $(V_p,\ldots,\psi_p)$ on $X$, where
$\psi'_p:(s_p\vert_{\pd V_p})^{-1}(0)\ra\pd X$ is given by
$\psi'_p:(q,B)\mapsto[\psi_p(q),(V_p,\ldots\psi_p),B]$ for
$(q,B)\in\pd V_p$ with $s_p(q)=0$. Then $\vdim(\pd X)=\vdim X-1$,
and $\pd X$ is compact if $X$ is compact.
\label{ug2def6}
\end{dfn}

In \S\ref{ug21} we explained that if $X$ is a manifold with
(g-)corners then there is a natural involution
$\si:\pd^2X\ra\pd^2X$. The same construction works for orbifolds,
and for Kuranishi spaces. That is, if $X$ is a Kuranishi space then
as in \cite[\S 2.6]{Joyc1} there is a natural strong diffeomorphism
$\bs\si:\pd^2X\ra\pd^2X$ with $\bs\si^2=\bs\id_X$. If $X$ is
oriented as in \S\ref{ug26} below then $\bs\si$ is
orientation-reversing.

\subsection{Fibre products of Kuranishi spaces}
\label{ug25}

We can define {\it fibre products\/} of Kuranishi spaces
\cite[Def.~2.28]{Joyc1}, as for fibre products of manifolds and
orbifolds in~\S\ref{ug21}.

\begin{dfn} Let $X,X'$ be Kuranishi spaces, $Y$ an orbifold, and
$\bs f:X\ra Y$, $\bs f':X'\ra Y$ be strongly smooth maps inducing
continuous maps $f:X\ra Y$ and $f':X'\ra Y$. Suppose at least one of
$\bs f,\bs f'$ is a strong submersion. We shall define the {\it
fibre product\/} $X\t_YX'$ or $X\t_{\bs f,Y,\bs f'}X'$, a Kuranishi
space. As a set, the underlying topological space $X\t_YX'$ is given
by~\eq{ug2eq3}.

Let $p\in X$, $p'\in X'$ and $q\in Y$ with $f(p)=q=f'(p')$. Let
$(V_p,E_p,s_p,\psi_p)$, $(V'_{\smash{p'}},E'_{\smash{p'}},
s'_{\smash{p'}},\psi'_{\smash{p'}})$ be sufficiently small Kuranishi
neighbourhoods in the germs at $p,p'$ in $X,X'$, and $f_p:V_p\ra Y$,
$f'_{\smash{p'}}:V'_{\smash{p'}}\ra Y$ be smooth maps in the germs
of $\bs f,\bs f'$ at $p,p'$ respectively. Define a Kuranishi
neighbourhood on $X\t_YX'$ by
\e
\begin{split}
\bigl(V_p\t_{f_p,Y,f'_{\smash{p'}}}&V'_{\smash{p'}},\pi_{V_p}^*(E_p)\op
\pi_{V'_{\smash{p'}}}^*(E'_{\smash{p'}}),\\
&(s_p\ci\pi_{V_p})\op(s'_{\smash{p'}}\ci\pi_{V'_{\smash{p'}}}),
(\psi_p\ci\pi_{V_p})\t(\psi'_{\smash{p'}}\ci\pi_{V'_{\smash{p'}}})\t
\chi_{pp'}\bigr).
\end{split}
\label{ug2eq5}
\e

Here $V_p\t_{f_p,Y,f'_{\smash{p'}}}V'_{\smash{p'}}$ is the fibre
product of orbifolds, and $\pi_{V_p},\pi_{V'_{\smash{p'}}}$ are the
projections from $V_p\t_YV'_{\smash{p'}}$ to $V_p,V'_{\smash{p'}}$.
The final term $\chi_{pp'}$ in \eq{ug2eq5} maps the biquotient terms
in \eq{ug2eq3} for $V_p\t_YV'_{\smash{p'}}$ to the same terms in
\eq{ug2eq3} for the set $X\t_YX'$. Coordinate changes between
Kuranishi neighbourhoods in $X,X'$ induce coordinate changes between
neighbourhoods \eq{ug2eq5}. So the systems of germs of Kuranishi
neighbourhoods and coordinate changes on $X,X'$ induce such systems
on $X\t_YX'$. This gives a {\it Kuranishi structure\/} on $X\t_YX'$,
making it into a {\it Kuranishi space}. Clearly
$\vdim(X\t_YX')=\vdim X+\vdim X'-\dim Y$, and $X\t_YX'$ is compact
if $X,X'$ are compact.
\label{ug2def7}
\end{dfn}

\subsection{Orientations and orientation conventions}
\label{ug26}

In \cite[\S 2.7]{Joyc1} we define {\it orientations\/} on Kuranishi
spaces. Our definition is basically equivalent to Fukaya et al.\
\cite[Def.~A1.17]{FOOO}, noting that our Kuranishi spaces correspond
to their Kuranishi spaces with a tangent bundle.

\begin{dfn} Let $X$ be a Kuranishi space. An {\it orientation\/} on $X$
assigns, for all $p\in X$ and all sufficiently small Kuranishi
neighbourhoods $(V_p,E_p,s_p,\psi_p)$ in the germ at $p$,
orientations on the fibres of the orbibundle $TV_p\op E_p$ varying
continuously over $V_p$. These must be compatible with coordinate
changes, in the following sense. Let $q\in\Im\psi_p$,
$(V_q,\ldots,\psi_q)$ be sufficiently small in the germ at $q$, and
$(\phi_{pq},\hat\phi_{pq})$ be the coordinate change from
$(V_q,\ldots,\psi_q)$ to $(V_p,\ldots,\psi_p)$ in the germ. Define
$\d\hat s_p$ near $s_q^{-1}(0)\subseteq V_q$ as in~\eq{ug2eq4}.

Locally on $V_q$, choose any orientation for the fibres of
$\phi_{pq}^*(TV_p)/(\d\phi_{pq})(TV_q)$, and let
$\phi_{pq}^*(E_p)/\hat\phi_{pq}(E_q)$ have the orientation induced
from this by the isomorphism $\d\hat s_p$ in \eq{ug2eq4}. These
induce an orientation on $\frac{\phi_{pq}^*(TV_p)}{
(\d\phi_{pq})(TV_q)}\op\frac{\phi_{pq}^*(E_p)}{
\hat\phi_{pq}(E_q)}$, which is independent of the choice for
$\phi_{pq}^*(TV_p)/(\d\phi_{pq})(TV_q)$. Thus, these local choices
induce a natural orientation on the orbibundle
$\frac{\phi_{pq}^*(TV_p)}{(\d\phi_{pq})(TV_q)}\op
\frac{\phi_{pq}^*(E_p)}{\hat\phi_{pq}(E_q)}$ near $s_q^{-1}(0)$. We
require that in oriented orbibundles over $V_q$ near $s_q^{-1}(0)$,
we have
\e
\begin{split}
\phi_{pq}^*\bigl[TV_p\op E_p\bigr]\cong (-1)^{\dim V_q(\dim
V_p-\dim V_q)}\bigl[TV_q\op E_q\bigr]\op&\\
\bigl[\ts\frac{\phi_{pq}^*(TV_p)}{(\d\phi_{pq})(TV_q)}\op
\frac{\phi_{pq}^*(E_p)}{\hat\phi_{pq}(E_q)}\bigr]&,
\end{split}
\label{ug2eq6}
\e
where $TV_p\op E_p$ and $TV_q\op E_q$ have the orientations assigned
by the orientation on $X$. An {\it oriented Kuranishi space\/} is a
Kuranishi space with an orientation.
\label{ug2def8}
\end{dfn}

Suppose $X,X'$ are oriented Kuranishi spaces, $Y$ is an oriented
orbifold, and $\bs f:X\ra Y$, $\bs f':X'\ra Y$ are strong
submersions. Then by \S\ref{ug24}--\S\ref{ug25} we have Kuranishi
spaces $\pd X$ and $X\t_YX'$. These can also be given orientations
in a natural way. We shall follow the orientation conventions of
Fukaya et al.~\cite[\S 45]{FOOO}.

\begin{conv} First, our conventions for manifolds:
\begin{itemize}
\setlength{\itemsep}{0pt}
\setlength{\parsep}{0pt}
\item[(a)] Let $X$ be an oriented manifold with boundary $\pd
X$. Then we define the orientation on $\pd X$ such that
$TX\vert_{\pd X}=\R_{\rm out}\op T(\pd X)$ is an isomorphism of
oriented vector spaces, where $\R_{\rm out}$ is oriented by an
outward-pointing normal vector to~$\pd X$.
\item[(b)] Let $X,X',Y$ be oriented manifolds, and $f:X\ra
Y$, $f':X'\ra Y$ be submersions. Then $\d f:TX\ra f^*(TY)$ and $\d
f':TX'\ra (f')^*(TY)$ are surjective maps of vector bundles over
$X,X'$. Choosing Riemannian metrics on $X,X'$ and identifying the
orthogonal complement of $\Ker\d f$ in $TX$ with the image $f^*(TY)$
of $\d f$, and similarly for $f'$, we have isomorphisms of vector
bundles over $X,X'$:
\e
TX\cong \Ker\d f\op f^*(TY) \quad\text{and}\quad TX'\cong
(f')^*(TY)\op\Ker\d f'.
\label{ug2eq7}
\e

Define orientations on the fibres of $\Ker\d f$, $\Ker\d f'$ over
$X,X'$ such that \eq{ug2eq7} are isomorphisms of oriented vector
bundles, where $TX,TX'$ are oriented by the orientations on $X,X'$,
and $f^*(TY),(f')^*(TY)$ by the orientation on $Y$. Then we define
the orientation on $X\t_YX'$ so that
\begin{equation*}
T(X\t_YX')\cong \pi_X^*(\Ker\d f)\op
(f\ci\pi_X)^*(TY)\op\pi_{X'}^*(\Ker\d f')
\end{equation*}
is an isomorphism of oriented vector bundles. Here $\pi_X:X\t_YX'\ra
X$ and $\pi_{X'}:X\t_YX'\ra X'$ are the natural projections, and
$f\ci\pi_X\equiv f'\ci\pi_{X'}$.
\end{itemize}
These extend immediately to orbifolds. They also extend to the
Kuranishi space versions in Definitions \ref{ug2def6} and
\ref{ug2def7}; for Definition \ref{ug2def7} they are described in
\cite[Conv.~45.1(4)]{FOOO}. An algorithm to deduce Kuranishi space
orientation conventions from manifold ones is described in~\cite[\S
2.7]{Joyc1}.
\label{ug2conv}
\end{conv}

If $X$ is an oriented Kuranishi space, we often write $-X$ for the
same Kuranishi space with the opposite orientation. Here is
\cite[Prop.~2.31]{Joyc1}, largely taken from Fukaya et
al.~\cite[Lem.~45.3]{FOOO}.

\begin{prop} Let\/ $X_1,X_2,\ldots$ be oriented Kuranishi spaces,
$Y,Y_1,\ldots$ be oriented orbifolds, and\/ $\bs f_1:X_1\ra
Y,\ldots$ be strongly smooth maps, with at least one strong
submersion in each fibre product below. Then the following hold, in
oriented Kuranishi spaces:
\begin{itemize}
\setlength{\itemsep}{0pt}
\setlength{\parsep}{0pt}
\item[{\rm(a)}] If\/ $\pd Y=\es,$ for\/ $\bs f_1:X_1\ra Y$ and\/
$\bs f_2:X_2\ra Y$ we have
\begin{equation}
\begin{gathered}
\pd(X_1\t_YX_2)=(\pd X_1)\t_YX_2\amalg (-1)^{\vdim X_1+\dim
Y}X_1\t_Y(\pd X_2)\\ \text{and}\quad X_1\t_YX_2=(-1)^{(\vdim
X_1-\dim Y)(\vdim X_2-\dim Y)}X_2\t_YX_1.
\end{gathered}
\label{ug2eq8}
\end{equation}
\item[{\rm(b)}] For\/ $\bs f_1:X_1\ra Y_1,$ $\bs f_2:X_2\ra Y_1\t
Y_2$ and\/ $\bs f_3:X_3\ra Y_2,$ we have
\begin{equation}
(X_1\t_{Y_1}X_2)\t_{Y_2}X_3=X_1\t_{Y_1}(X_2\t_{Y_2}X_3).
\label{ug2eq9}
\end{equation}
\item[{\rm(c)}] For\/ $\bs f_1:X_1\ra Y_1\t Y_2,$ $\bs f_2:X_2\ra
Y_1$ and\/ $\bs f_3:X_3\ra Y_2,$ we have
\end{itemize}
\begin{equation}
X_1\t_{Y_1\t Y_2}(X_2\t X_3)=(-1)^{\dim Y_2(\dim Y_1+\vdim X_2)}
(X_1\t_{Y_1}X_2)\t_{Y_2}X_3.
\label{ug2eq10}
\end{equation}
\label{ug2prop}
\end{prop}

\subsection{Coorientations}
\label{ug27}

To define Kuranishi cohomology in \S\ref{ug4} we we will also need a
notion of {\it relative orientation} for (strong) submersions. We
call this a {\it coorientation},~\cite[\S 2.8]{Joyc1}.

\begin{dfn} Let $X,Y$ be orbifolds, and $f:X\ra Y$ a submersion. A
{\it coorientation\/} for $(X,f)$ is a choice of orientations on the
fibres of the vector bundle $\Ker\d f$ over $X$ which vary
continuously over $X$. Here $\d f:TX\ra f^*(TY)$ is the derivative
of $f$, a morphism of vector bundles, which is surjective as $f$ is
a submersion. Thus $\Ker\d f$ is a vector bundle over $X$, of
rank~$\dim X-\dim Y$.

Now let $X$ be a Kuranishi space, $Y$ an orbifold, and $\bs f:X\ra
Y$ a strong submersion. A {\it coorientation\/} for $(X,\bs f)$
assigns, for all $p\in X$ and all sufficiently small Kuranishi
neighbourhoods $(V_p,E_p,s_p,\psi_p)$ in the germ at $p$ with
submersion $f_p:V_p\ra Y$ representing $\bs f$, orientations on the
fibres of the orbibundle $\Ker\d f_p\op E_p$ varying continuously
over $V_p$, where $\d f_p:TV_p\ra f_p^*(TY)$ is the (surjective)
derivative of $f_p$.

These must be compatible with coordinate changes, in the following
sense. Let $q\in\Im\psi_p$, $(V_q,\ldots,\psi_q)$ be sufficiently
small in the germ at $q$, let $f_q:V_q\ra Y$ represent $\bs f$, and
$(\phi_{pq},\hat\phi_{pq})$ be the coordinate change from
$(V_q,\ldots,\psi_q)$ to $(V_p,\ldots,\psi_p)$ in the germ. Then we
require that in oriented orbibundles over $V_q$ near $s_q^{-1}(0)$,
we have
\e
\begin{split}
\phi_{pq}^*\bigl[\Ker\d f_p\op E_p\bigr]\cong (-1)^{\dim V_q(\dim
V_p-\dim V_q)}\bigl[\Ker\d f_q\op E_q\bigr]\op&\\
\bigl[\ts\frac{\phi_{pq}^*(TV_p)}{(\d\phi_{pq})(TV_q)}\op
\frac{\phi_{pq}^*(E_p)}{\hat\phi_{pq}(E_q)}\bigr]&,
\end{split}
\label{ug2eq11}
\e
by analogy with \eq{ug2eq6}, where $\frac{\phi_{pq}^*(TV_p)}{
(\d\phi_{pq})(TV_q)}\op\frac{\phi_{pq}^*(E_p)}{\hat\phi_{pq}(E_q)}$
is oriented as in Definition~\ref{ug2def8}.
\label{ug2def9}
\end{dfn}

Suppose now that $Y$ is oriented. Then an orientation on $X$ is
equivalent to a coorientation for $(X,\bs f)$, since for all $p,
(V_p,\ldots,\psi_p),f_p$ as above, the isomorphism $TV_p\cong
f_p^*(TY)\op \Ker\d f_p$ induces isomorphisms of orbibundles
over~$V_p$:
\e
\bigl(TV_p\op E_p\bigr)\cong f_p^*(TY)\op\bigl(\Ker\d f_p\op
E_p\bigr).
\label{ug2eq12}
\e
There is a 1-1 correspondence between orientations on $X$ and
coorientations for $(X,\bs f)$ such that \eq{ug2eq12} holds in
oriented orbibundles, where $TV_p\op E_p$ is oriented by the
orientation on $X$, and $\Ker\d f_p\op E_p$ by the coorientation for
$(X,\bs f)$, and $f_p^*(TY)$ by the orientation on $Y$. Taking the
direct sum of $f_p^*(TY)$ with \eq{ug2eq11} and using \eq{ug2eq12}
yields \eq{ug2eq6}, so this is compatible with coordinate changes.

In \cite[Conv.~2.33]{Joyc1} we give our conventions for
coorientations of boundaries and fibre products. These correspond to
Convention \ref{ug2conv} under the 1-1 correspondence between
orientations on $X$ and coorientations for $(X,\bs f)$ above when
$Y$ is oriented. So the analogue of Proposition \ref{ug2prop} holds
for coorientations. In particular, for strong submersions $\bs
f_a:X_a\ra Y$ with $(X_a,\bs f_a)$ cooriented for $a=1,2,3$, taking
$\pd Y=\es$ in \eq{ug2eq13}, we have
\begin{gather}
\begin{split}
\bigl(\pd(X_1\t_YX_2),\bs\pi_Y\bigr)&\cong\bigl((\pd
X_1)\t_YX_2,\bs\pi_Y\bigr)\,\amalg \\
&\qquad\quad (-1)^{\vdim X_1+\dim Y} \bigl(X_1\t_Y(\pd
X_2),\bs\pi_Y\bigr),
\end{split}
\label{ug2eq13}\\
\bigl(X_1\t_YX_2,\bs\pi_Y\bigr)\cong(-1)^{(\vdim X_1-\dim Y)(\vdim
X_2-\dim Y)}\bigl(X_2\t_YX_1,\bs\pi_Y\bigr),
\label{ug2eq14}\\
\bigl((X_1\t_YX_2)\t_YX_3,\bs\pi_Y\bigr)\cong
\bigl(X_1\t_Y(X_2\t_YX_3),\bs\pi_Y\bigr).
\label{ug2eq15}
\end{gather}
Similarly, if $X_1$ is oriented, $\bs f_1:X_1\ra Y$ is strongly
smooth, $\bs f_2:X_2\ra Y$ is a cooriented strong submersion, and
$\pd Y=\es$ then
\e
\pd(X_1\t_YX_2)\cong \bigl((\pd X_1)\t_YX_2\bigr)\amalg (-1)^{\vdim
X_1+\dim Y}\bigl(X_1\t_Y(\pd X_2)\bigr)
\label{ug2eq16}
\e
in oriented Kuranishi spaces.

\section{Kuranishi homology}
\label{ug3}

{\it Kuranishi homology\/} \cite[\S 4]{Joyc1} is a homology theory
of orbifolds $Y$ in which the chains are isomorphism classes $[X,\bs
f,\bs G]$, where $X$ is a compact, oriented Kuranishi space, $\bs
f:X\ra Y$ is strongly smooth, and $\bs G$ is some extra data called
{\it gauge-fixing data}. It is isomorphic to singular homology
$H_*^\rsi(Y;R)$.

\subsection{Gauge-fixing data}
\label{ug31}

Let $X$ be a compact Kuranishi space, $Y$ an orbifold, and $\bs
f:X\ra Y$ a strongly smooth map. A key ingredient in the definition
of Kuranishi homology in \cite{Joyc1} is the idea of {\it gauge
fixing data\/} $\bs G$ for $(X,\bs f)$ studied in \cite[\S
6]{Joyc1}. Define $P=\coprod_{k=0}^\iy \R^k/S_k$, where the
symmetric group $S_k$ acts on $\R^k$ by permuting the coordinates
$x_1,\ldots,x_k$. For $n=0,1,2,\ldots$, define $P_n\subset P$ by
$P_n=\coprod_{k=0}^n\R^k/S_k$. Gauge-fixing data $\bs G$ for $(X,\bs
f)$ consists of a cover of $X$ by Kuranishi neighbourhoods
$(V^i,E^i,s^i,\psi^i)$ for $i$ in a finite indexing set $I$,
together with smooth maps $f^i:V^i\ra Y$ representing $\bs f$ and
maps $G^i:E^i\ra P_n\subset P$ for some $n\gg 0$, and continuous
partitions of unity $\eta_i:X\ra[0,1]$ and $\eta_i^j:V^j\ra[0,1]$,
satisfying many conditions. One important condition, responsible for
Theorem \ref{ug3thm1}(b) below, is that each $G^i:E^i\ra P$ should
be a {\it finite\/} map, that is, $(G^i)^{-1}(p)$ is finitely many
points for all~$p\in P$.

Users of Kuranishi homology do not need to know exactly what
gauge-fixing data is, so we will not define it. Here are the
important properties of gauge-fixing data, which are proved
in~\cite[\S 3]{Joyc1}.

\begin{thm} Consider pairs\/ $(X,\bs f)$, where\/ $X$ is a compact
Kuranishi space, $Y$ an orbifold, and\/ $\bs f:X\ra Y$ a strongly
smooth map. In\/ {\rm\cite[\S 3.1]{Joyc1}} we define
\begin{bfseries}gauge-fixing data\end{bfseries}\/ $\bs G$ for such
pairs $(X,\bs f),$ with the following properties:
\begin{itemize}
\setlength{\itemsep}{0pt}
\setlength{\parsep}{0pt}
\item[{\rm(a)}] Every pair $(X,\bs f)$ admits a (nonunique) choice
of gauge-fixing data $\bs G$. If\/ $\Ga\subseteq\Aut(X,\bs f)$ is a
finite subgroup then we can choose $\bs G$ to be $\Ga$-invariant.
\item[{\rm(b)}] For all pairs\/ $(X,\bs f)$ and choices of
gauge-fixing data\/ $\bs G$ for\/ $(X,\bs f),$ the automorphism
group\/ $\Aut(X,\bs f,\bs G)$ of isomorphisms\/ $(\bs a,\bs
b):(X,\bs f,\bs G)\ra(X,\bs f,\bs G)$ is finite.
\item[{\rm(c)}] Suppose $\bs G$ is gauge-fixing data for $(X,\bs
f)$ and\/ $\Ga$ is a finite group acting on $(X,\bs f,\bs G)$ by
isomorphisms. Then we can form the quotient\/ $\ti X=X/\Ga,$ a
compact Kuranishi space, with projection $\bs\pi:X\ra\ti X,$ and\/
$\bs f$ pushes down to $\bs{\ti f}:\ti X\ra Y$ with\/ $\bs f=\bs{\ti
f}\ci\bs\pi$. As in\/ {\rm\cite[\S 3.4]{Joyc1},} we can define
gauge-fixing data\/ $\bs{\ti G}$ for\/ $(\ti X,\bs{\ti f}),$ which
is the natural push down\/ $\bs\pi_*(\bs G)$ of\/ $\bs G$ to~$\ti
X$.
\item[{\rm(d)}] If\/ $\bs G$ is gauge-fixing data for\/ $(X,\bs
f),$ it has a restriction\/ $\bs G\vert_{\pd X}$ defined in\/
{\rm\cite[\S 3.5]{Joyc1},} which is gauge-fixing data for\/~$(\pd
X,\bs f\vert_{\pd X})$.
\item[{\rm(e)}] Let\/ $X$ be a compact, oriented Kuranishi space
with corners (not g-corners), $\bs f:X\ra Y$ be strongly smooth,
and\/ $\bs\si:\pd^2X\ra\pd^2X$ be the natural involution described
in\/ {\rm\S\ref{ug24}}. Suppose $\bs H$ is gauge-fixing data for
$(\pd X,\bs f\vert_{\pd X})$. Then there exists gauge-fixing data\/
$\bs G$ for\/ $(X,\bs f)$ with\/ $\bs G\vert_{\pd X}=\bs H$ if and
only if\/ $\bs H\vert_{\pd^2X}$ is invariant under\/ $\bs\si$. If
also\/ $\Ga$ is a finite subgroup of\/ $\Aut(X,\bs f),$ and\/ $\bs
H$ is invariant under\/ $\Ga\vert_{\pd X},$ then we can choose\/
$\bs G$ to be\/ $\Ga$-invariant.
\item[{\rm(f)}] Let\/ $Y,Z$ be orbifolds, and\/ $h:Y\ra Z$ a smooth
map. Suppose $X$ is a compact Kuranishi space, $\bs f:X\ra Y$ is
strongly smooth, and\/ $\bs G$ is gauge-fixing data for $(X,\bs f)$.
Then as in {\rm\cite[\S 3.7]{Joyc1},} we can define gauge-fixing
data\/ $h_*(\bs G)$ for $(X,h\ci\bs f)$. It satisfies $(g\ci
h)_*(\bs G)=g_*(h_*(\bs G))$.
\end{itemize}
\label{ug3thm1}
\end{thm}

In order to define a working homology theory, perhaps the most
important property is Theorem \ref{ug3thm1}(b). In \cite[\S
4.9]{Joyc1} we show that if we define {\it na\"\i ve Kuranishi
(co)homology} $KH_*^\na,KH^*_\na(Y;R)$ as in \S\ref{ug32} and
\S\ref{ug42} but omitting all (co-)gauge-fixing data, then
$KH_*^\na(Y;R)=0=KH^*_\na(Y;R)$ for all orbifolds $Y$ and
$\Q$-algebras $R$. For compact $Y$, we do this by constructing an
explicit cochain whose boundary is the identity cocycle $[Y,\id_Y]$
in $KC^*_\na(Y;R)$. This explicit cochain involves a cocycle $[X,\bs
f]$ whose automorphism group $\Aut(X,\bs f)$ is infinite. Including
(co-)gauge-fixing data prevents this from happening, as it ensures
that all automorphism groups $\Aut(X,\bs f,\bs G)$ are finite.

\subsection{Kuranishi homology}
\label{ug32}

We can now define the {\it Kuranishi homology\/} of an orbifold,
\cite[\S 4.2]{Joyc1}.

\begin{dfn} Let $Y$ be an orbifold. Consider triples $(X,\bs f,\bs
G)$, where $X$ is a compact, oriented Kuranishi space, $\bs f:X\ra
Y$ is strongly smooth, and $\bs G$ is gauge-fixing data for $(X,\bs
f)$. Write $[X,\bs f,\bs G]$ for the isomorphism class of $(X,\bs
f,\bs G)$ under isomorphisms $(\bs a,\bs b):(X,\bs f,\bs G)\ra(\ti
X,\bs{\ti f},\bs{\ti G})$, where $\bs a$ must identify the
orientations of $X,\ti X$, and $\bs b$ lifts $\bs a$ to the
Kuranishi neighbourhoods $(V^i,\ldots,\psi^i),(\ti
V^i,\ldots,\ti\psi^i)$ in~$\bs G,\bs{\ti G}$.

Let $R$ be a $\Q$-algebra, for instance $\Q,\R$ or $\C$. For each
$k\in\Z$, define $KC_k(Y;R)$ to be the $R$-module of finite
$R$-linear combinations of isomorphism classes $[X,\bs f,\bs G]$ for
which $\vdim X=k$, with the relations:
\begin{itemize}
\setlength{\itemsep}{0pt}
\setlength{\parsep}{0pt}
\item[(i)] Let $[X,\bs f,\bs G]$ be an isomorphism class,
and write $-X$ for $X$ with the opposite orientation. Then in
$KC_k(Y;R)$ we have
\begin{equation*}
[X,\bs f,\bs G]+[-X,\bs f,\bs G]=0.
\end{equation*}
\item[(ii)] Let $[X,\bs f,\bs G]$ be an isomorphism class. Suppose
that $X$ may be written as a disjoint union $X=X_+\amalg X_-$ of
compact, oriented Kuranishi spaces, and that for each Kuranishi
neighbourhood $(V^i,\ldots,\psi^i)$ for $i\in I$ in $\bs G$ we may
write $V^i=V^i_+\amalg V^i_-$ for open and closed subsets $V^i_\pm$
of $V^i$, such that $\Im\psi^i\vert_{V^i_+}\subseteq X_+$ and
$\Im\psi^i\vert_{V^i_-}\subseteq X_-$. Then we may define
gauge-fixing data $\bs G\vert_{X_\pm}$ for $(X_\pm,\bs
f\vert_{X_\pm})$, with Kuranishi neighbourhoods $(V^i_\pm,
E^i\vert_{V^i_\pm},s^i\vert_{V^i_\pm},\psi^i\vert_{V^i_\pm})$ for
$i\in I$ with $V^i_\pm\ne\emptyset$. In $KC_k(Y;R)$ we have
\begin{equation*}
[X,\bs f,\bs G]=[X_+,\bs f\vert_{X_+},\bs G\vert_{X_+}]+[X_-,\bs
f\vert_{X_-},\bs G\vert_{X_-}].
\end{equation*}
\item[(iii)] Let $[X,\bs f,\bs G]$ be an isomorphism class, and
suppose $\Ga$ is a finite group acting on $(X,\bs f,\bs G)$ by
orientation-preserving automorphisms. Then $\ti X=X/\Ga$ is a
compact, oriented Kuranishi space, with a projection $\bs\pi:X\ra\ti
X$. As in Theorem \ref{ug3thm1}(c), $\bs f,\bs G$ push down to a
strong submersion $\bs\pi_*(\bs f)=\bs{\ti f}:\ti X\ra Y$ and
gauge-fixing data $\bs\pi_*(\bs G)=\bs{\ti G}$ for $(\ti X,\bs{\ti
f})$. Then
\begin{equation*}
\bigl[X/\Ga,\bs\pi_*(\bs f),\bs\pi_*(\bs G)\bigr]=\frac{1}{\md{\Ga}}
\,\bigl[X,\bs f,\bs G\bigr]
\end{equation*}
\end{itemize}
in $KC_k(Y;R)$. Elements of $KC_k(Y;R)$ will be called {\it
Kuranishi chains}.

Define the {\it boundary operator\/} $\pd:KC_k(Y;R)\ra
KC_{k-1}(Y;R)$ by
\begin{equation*}
\pd:\ts\sum_{a\in A}\rho_a[X_a,\bs f_a,\bs G_a]\longmapsto
\ts\sum_{a\in A}\rho_a[\pd X_a,\bs f_a\vert_{\pd X_a},\bs G_a
\vert_{\pd X_a}],
\end{equation*}
where $A$ is a finite indexing set and $\rho_a\in R$ for $a\in A$.
This is a morphism of $R$-modules. Clearly, $\pd$ takes each
relation (i)--(iii) in $KC_k(Y;R)$ to the corresponding relation in
$KC_{k-1}(Y;R)$, and so $\pd$ is well-defined.

Recall from \S\ref{ug21} and \S\ref{ug24} that if $X$ is an oriented
Kuranishi space then there is a natural orientation-reversing strong
diffeomorphism $\bs\si:\pd^2X\ra\pd^2X$, with $\bs\si^2=
\bs\id_{\pd^2X}$. If $[X,\bs f,\bs G]$ is an isomorphism class then
this $\bs\si$ extends to an isomorphism $(\bs\si,\bs\tau)$ of
$(\pd^2X,\bs f\vert_{\pd^2X},\bs G\vert_{\pd^2X})$. So part (i) in
$KC_{k-2}(Y;R)$ yields
\e
[\pd^2X,\bs f\vert_{\pd^2X},\bs G\vert_{\pd^2X}]+[\pd^2X,\bs
f\vert_{\pd^2X},\bs G\vert_{\pd^2X}]=0 \quad\text{in
$KC_{k-2}(Y;R)$.}
\label{ug3eq1}
\e
As $R$ is a $\Q$-algebra we may multiply \eq{ug3eq1} by $\ha$ to get
$[\pd^2X,\bs f\vert_{\pd^2X},\bs G\vert_{\pd^2X}]=0$. Therefore
$\pd\ci\pd=0$ as a map~$KC_k(Y;R)\ra KC_{k-2}(Y;R)$.

Define the {\it Kuranishi homology group\/} $KH_k(Y;R)$ of $Y$ for
$k\in\Z$ to be
\begin{equation*}
KH_k(Y;R)=\frac{\Ker\bigl(\pd:KC_k(Y;R)\ra KC_{k-1}(Y;R)\bigr)}{
\Im\bigl(\pd:KC_{k+1}(Y;R)\ra KC_k(Y;R)\bigr)}\,.
\end{equation*}

Let $Y,Z$ be orbifolds, and $h:Y\ra Z$ a smooth map. Define the {\it
pushforward} $h_*:KC_k(Y;R)\ra KC_k(Z;R)$ for $k\in\Z$ by
\begin{equation*}
h_*:\ts\sum_{a\in A}\rho_a\bigl[X_a,\bs f_a,\bs G_a\bigr]\longmapsto
\ts\sum_{a\in A}\rho_a\bigl[X_a,h\ci\bs f_a,h_*(\bs G_a)\bigr],
\end{equation*}
with $h_*(\bs G_a)$ as in Theorem \ref{ug3thm1}(f). These take
relations (i)--(iii) in $KC_k(Y;R)$ to (i)--(iii) in $KC_k(Z;R)$,
and so are well-defined. They satisfy $h_*\ci\pd=\pd\ci h_*$, so
they induce morphisms of homology groups $h_*:KH_k(Y;R)\ra
KH_k(Z;R)$. Pushforward is functorial, that is, $(g\ci h)_*=g_*\ci
h_*$, on chains and homology.
\label{ug3def}
\end{dfn}

\subsection{Singular homology and Kuranishi homology}
\label{ug33}

Let $Y$ be an orbifold, and $R$ a $\Q$-algebra. Then we can define
the {\it singular homology groups\/} $H_k^\rsi(Y;R)$, as in Bredon
\cite[\S IV]{Bred}. Write $C_k^\rsi(Y;R)$ for the $R$-module spanned
by {\it smooth singular $k$-simplices\/} in $Y$, which are smooth
maps $\si:\De_k\ra Y$, where $\De_k$ is the $k$-simplex
\begin{equation*}
\De_k=\bigl\{(x_0,\ldots,x_k)\in\R^{k+1}:x_i\ge 0,\;\>
x_0+\ldots+x_k=1\bigr\}.
\end{equation*}
As in \cite[\S IV.1]{Bred}, the boundary operator
$\pd:C_k^\rsi(Y;R)\ra C_{k-1}^\rsi(Y;R)$ is given by
\begin{equation*}
\pd:\ts\sum_{a\in A}\rho_a\,\si_a\longmapsto \ts\sum_{a\in
A}\sum_{j=0}^k(-1)^j\rho_a(\si_a\ci F_j^k),
\end{equation*}
where $F_j^k:\De_{k-1}\ra\De_k$, $F_j^k:(x_0,\ldots,x_{k-1})
\mapsto(x_0,\ldots,x_{j-1},0,x_j,\ldots,x_{k-1})$ for
$j=0,\ldots,k$. Then $\pd^2=0$, and $H_*^\rsi(Y;R)$ is the homology
of~$\bigl(C_k^\rsi(Y;R),\pd\bigr)$.

Define $R$-module morphisms $C_k^\rsi(Y;R)\ra KC_k(Y;R)$ for $k\ge
0$ by
\begin{equation*}
\ts\Pi_\rsi^\Kh:\sum_{a\in A}\rho_a\si_a\longmapsto \ts\sum_{a\in
A}\rho_a\bigl[\De_k,\si_a,\bs G_{\De_k}\bigr],
\end{equation*}
where $\bs G_{\De_k}$ is an explicit choice of gauge-fixing data for
$(\De_k,\si_a)$ given in \cite[\S 4.3]{Joyc1}. These satisfy
$\pd\ci\Pi_\rsi^\Kh=\Pi_\rsi^\Kh\ci\pd$, and so induce $R$-module
morphisms
\e
\Pi_\rsi^\Kh:H_k^\rsi(Y;R)\longra KH_k(Y;R).
\label{ug3eq2}
\e
Here is \cite[Cor.~4.10]{Joyc1}, one of the main results
of~\cite{Joyc1}.

\begin{thm} Let\/ $Y$ be an orbifold and\/ $R$ a $\Q$-algebra. Then
$\Pi_\rsi^\Kh$ in\/ \eq{ug3eq2} is an isomorphism, so that\/
$KH_k(Y;R)\cong H_k^\rsi(Y;R),$ with\/ $KH_k(Y;R)=\{0\}$
when\/~$k<0$.
\label{ug3thm2}
\end{thm}

The proof of Theorem \ref{ug3thm2} in \cite[App.~A--C]{Joyc1} is
very long and complex, taking up a third of \cite{Joyc1}. The
problem is to construct an inverse for $\Pi_\rsi^\Kh$ in
\eq{ug3eq2}. This is related to Fukaya and Ono's construction of
{\it virtual cycles\/} for compact, oriented Kuranishi spaces
without boundary in \cite[\S 6]{FuOn}, and uses some of the same
ideas. But dealing with boundaries and corners of the Kuranishi
spaces in Kuranishi chains, and the relations in the Kuranishi chain
groups $KC_*(Y;R)$, increases the complexity by an order of
magnitude.

The basic idea of the proof is to take a class $\al\in KH_k(Y;R)$
and represent it by cycles $\sum_{a\in A}\rho_a[X_a,\bs f_a,\bs
G_a]$ with better and better properties, until eventually we
represent it by a cycle in the image of $\Pi_\rsi^\Kh:
C_k^\rsi(Y;R)\ra KC_k(Y;R)$, so showing that \eq{ug3eq2} is
surjective. Here (somewhat oversimplified) are the main steps:
firstly, by `cutting' the $X_a$ into small pieces $X_{ac}$ for $c\in
C_a$, we show we can represent $\al$ by a cycle $\sum_{a\in
A}\sum_{c\in C_a}\rho_a[X_{ac},\bs f_{ac},\bs G_{ac}]$ such that
$(X_{ac},\bs f_{ac},\bs G_{ac})$ is the quotient of a triple
$(\acute X_{ac},\bs{\acute f}_{ac},\bs{ \acute G}_{ac})$ by a finite
group $\Ga_{ac}$, where $\acute X_{ac}$ has {\it trivial
stabilizers}.

Thus by Definition \ref{ug3def}(iii) we can represent $\al$ by a
cycle $\sum_{a,c}\rho_a\md{\Ga_{ac}}^{-1}[\acute X_{ac},\ab
\smash{\bs{\acute f}_{ac},\bs{\acute G}_{ac}}]$ involving only
Kuranishi spaces $\smash{\acute X_{ac}}$ with trivial stabilizers.
Such spaces can be deformed to manifolds with g-corners (by
single-valued perturbations, not multisections). So secondly, we
show we can represent $\al$ by a cycle $\sum_{a,c}\rho_a
\md{\Ga_{ac}}^{-1}[\ti X_{ac},\ti f_{ac},\bs{\ti G}_{ac}]$ in which
the $\ti X_{ac}$ are manifolds, and $\smash{\ti f_{ac}:\ti X_{ac}\ra
Y}$ are smooth maps. Then thirdly, we triangulate the $\smash{\ti
X_{ac}}$ by simplices $\De_k$, and so prove that we can represent
$\al$ by a cycle in the image of $\Pi_\rsi^\Kh$. In this third step
it is vital to work with manifolds {\it with g-corners}, as in
\S\ref{ug21}, not just manifolds with corners, since otherwise we
would not be able to construct the homology between
$\sum_{a,c}\rho_a\md{\Ga_{ac}}^{-1} [\ti X_{ac},\ti f_{ac},\bs{\ti
G}_{ac}]$ and the singular cycle.

In the proof we use the fact that $R$ is a $\Q$-{\it algebra\/} in
two different ways. When we replace $[X_{ac},\bs f_{ac},\bs G_{ac}]$
by $\md{\Ga_{ac}}^{-1}[\acute X_{ac},\bs{\acute f}_{ac},\bs{\acute
G}_{ac}]$ we must have $\md{\Ga_{ac}}^{-1}\in R$, so we need
$\Q\subseteq R$. And when we deform $\acute X_{ac}$ to manifolds
$\ti X_{ac}$, to make $\sum_{a,c}\rho_a\md{\Ga_{ac}}^{-1}[\ti
X_{ac},\ti f_{ac},\bs{\ti G}_{ac}]$ a cycle, we should ensure our
perturbations are preserved by the automorphism groups
$\smash{\Aut(\acute X_{ac},\bs{\acute f}_{ac},\bs{\acute G}_{ac})}$.
In fact this may not be possible, if $\Aut(\acute X_{ac},\bs{\acute
f}_{ac},\bs{\acute G}_{ac})$ has fixed points. So instead, we choose
one perturbation $\ti X_{ac}$, and then take the average of the
images of this perturbation under $\Aut(\acute X_{ac},\bs{\acute
f}_{ac},\bs{\acute G}_{ac})$. This requires us to divide by
$\bmd{\Aut(\acute X_{ac},\bs{\acute f}_{ac},\bs{\acute G}_{ac})}$,
so again we need $\Q\subseteq R$. Also, for this step it is
necessary that automorphism groups $\Aut(X,\bs f,\bs G)$ should be
{\it finite}, as in Theorem \ref{ug3thm1}(b), and this was the
reason for introducing gauge-fixing data.

The theorem means that in many problems, particularly areas of in
Symplectic Geometry involving moduli spaces of $J$-holomorphic
curves, we can use Kuranishi chains and homology instead of singular
chains and homology, which can simplify proofs considerably, and
also improve results.

\section{Kuranishi cohomology}
\label{ug4}

We now discuss the Poincar\'e dual theory of {\it Kuranishi
cohomology\/} $KH^*(Y;R)$, which is isomorphic to
compactly-supported cohomology $H^*_\cs(Y;R)$. It is defined using a
complex of Kuranishi cochains $KC^*(Y;R)$ spanned by isomorphism
classes $[X,\bs f,\bs C]$ of triples $(X,\bs f,\bs C)$, where $X$ is
a compact Kuranishi space, $\bs f:X\ra Y$ is a {\it cooriented
strong submersion}, and $\bs C$ is {\it co-gauge-fixing data}.

As is usual for cohomology, Kuranishi cohomology has an associative,
supercommutative {\it cup product\/} $\cup:KH^k(Y;R)\t KH^l(Y;R)\ra
KH^{k+l}(Y;R)$, and there is also a {\it cap product\/}
$\cap:KH_k(Y;R)\t KH^l(Y;R)\ra KH_{k-l}(Y;R)$ relating Kuranishi
homology and Kuranishi (co)homology, which makes $KH_*(Y;R)$ into a
module over $KH^*(Y;R)$. More unusually, we can define $\cup,\cap$
naturally on Kuranishi (co)chains, and $\cup$ is associative and
supercommutative on $KC^*(Y;R)$, and $\cap$ makes $KC_*(Y;R)$ into a
module over~$KC^*(Y;R)$.

\subsection{Co-gauge-fixing data}
\label{ug41}

Let $X$ be a compact Kuranishi space, $Y$ an orbifold, and $\bs
f:X\ra Y$ a strong submersion. Kuranishi cohomology is based on the
idea of {\it co-gauge-fixing data} $\bs C$ for $(X,\bs f)$. This is
very similar to gauge-fixing data $\bs G$ in \S\ref{ug31}, and
consists of a finite cover of $X$ by Kuranishi neighbourhoods
$(V^i,E^i,s^i,\psi^i)$ for $i\in I$, submersions $f^i:V^i\ra Y$
representing $\bs f$ and maps $C^i:E^i\ra P_n\subset P$ for some
$n\gg 0$, and partitions of unity $\eta_i:X\ra[0,1]$
and~$\eta_i^j:V^j\ra[0,1]$.

Here are the important properties of co-gauge-fixing data, which are
proved in \cite[\S 3]{Joyc1}. Part (g) makes cup products work on
Kuranishi cochains, and part (h) makes cap products work. It was
difficult to find a definition of (co-)gauge-fixing data for which
properties (a)--(h) all hold at once; a large part of the complexity
of \cite[\S 3]{Joyc1} is due to the author's determination to ensure
that cup products should be associative and supercommutative {\it at
the cochain level}. This is not essential for a well-behaved
(co)homology theory, but is extremely useful in the
applications~\cite{AkJo,Joyc2,Joyc3}.

\begin{thm} Consider pairs\/ $(X,\bs f),$ where\/ $X$ is a compact
Kuranishi space, $Y$ an orbifold, and\/ $\bs f:X\ra Y$ a strong
submersion. In\/ {\rm\cite[\S 3.1]{Joyc1}} we define
\begin{bfseries}co-gauge-fixing data\end{bfseries}\/ $\bs C$ for such
pairs $(X,\bs f)$. It satisfies the analogues of Theorem
{\rm\ref{ug3thm1}(a)--(e),} and also:
\begin{itemize}
\setlength{\itemsep}{0pt}
\setlength{\parsep}{0pt}
\item[{\rm(f)}] Let\/ $Y,Z$ be orbifolds, and\/ $h:Y\ra Z$ a smooth,
proper map. Suppose $X$ is a compact Kuranishi space, $\bs f:X\ra Z$
is a strong submersion, and\/ $\bs C$ is co-gauge-fixing data for
$(X,\bs f)$. Then the fibre product\/ $Y\t_{h,Z,\bs f}X$ is a
compact Kuranishi space, and\/ $\bs\pi_Y:Y\t_ZX\ra Y$ is a strong
submersion. As in {\rm\cite[\S 3.7]{Joyc1},} we can define
co-gauge-fixing data\/ $h^*(\bs C)$ for $(Y\t_ZX,\bs\pi_Y)$. It
satisfies $(g\ci h)^*(\bs C)=h^*(g^*(\bs C))$.
\item[{\rm(g)}] Let\/ $X_1,X_2,X_3$ be compact Kuranishi spaces, $Y$
an orbifold, $\bs f_i:X_i\ra Y$ be strong submersions for\/
$i=1,2,3,$ and\/ $\bs C_i$ be co-gauge-fixing data for\/ $(X_i,\bs
f_i)$ for\/ $i=1,2,3$. Then\/ {\rm\cite[\S 3.8]{Joyc1}} defines
co-gauge-fixing data $\bs C_1\t_Y\bs C_2$ for\/ $(X_1\t_{\bs
f_1,Y,\bs f_2}X_2,\bs\pi_Y)$ from\/~$\bs C_1,\bs C_2$.

This construction is \begin{bfseries}symmetric\end{bfseries}, in
that it yields isomorphic co-gauge-fixing data for\/
$(X_1\t_YX_2,\bs\pi_Y)$ and\/ $(X_2\t_YX_1,\bs\pi_Y)$ under the
natural isomorphism\/ $X_1\t_YX_2\cong X_2\t_YX_1$. It is also
\begin{bfseries}associative\end{bfseries}, in that it yields
isomorphic co-gauge-fixing data for\/
$\bigl((X_1\!\t_Y\!X_2)\!\t_Y\!X_3,\bs\pi_Y\bigr)$ and\/
$\bigl(X_1\!\t_Y\!(X_2\!\t_Y\!X_3),\bs\pi_Y \bigr)$
under\/~$(X_1\!\t_Y\!X_2)\!\t_Y\!X_3\!\cong\!X_1\!\t_Y\!
(X_2\!\t_Y\!X_3)$.

These properties also have straightforward generalizations to
multiple fibre products involving more than one orbifold\/ $Y,$ such
as\/ \eq{ug2eq9} and\/~\eq{ug2eq10}.
\item[{\rm(h)}] Let\/ $X_1,X_2$ be compact Kuranishi spaces, $Y$
an orbifold, $\bs f_1:X_1\ra Y$ be strongly smooth, $\bs f_2:X_2\ra
Y$ be a strong submersion, $\bs G_1$ be gauge-fixing data for
$(X_1,\bs f_1),$ and\/ $\bs C_2$ be co-gauge-fixing data for
$(X_2,\bs f_2)$. Then\/ {\rm\cite[\S 3.8]{Joyc1}} defines
gauge-fixing data $\bs G_1\t_Y\bs C_2$ for\/ $(X_1\t_{\bs f_1,Y,\bs
f_2}X_2,\bs\pi_Y)$ from\/ $\bs G_1,\bs C_2$. If also $\bs f_3:X_3\ra
Y$ is a strong submersion and\/ $\bs C_3$ is co-gauge-fixing data
for $(X_3,\bs f_3)$ then the natural isomorphism
$(X_1\!\t_Y\!X_2)\!\t_Y\!X_3\cong X_1\!\t_Y\!(X_2\!\t_Y\!X_3)$
identifies $(\bs G_1\!\t_Y\!\bs C_2)\!\t_Y\!\bs C_3$ and\/~$\bs
G_1\!\t_Y\!(\bs C_2\!\t_Y\!\bs C_3)$.
\end{itemize}
\label{ug4thm1}
\end{thm}

\subsection{Kuranishi cohomology}
\label{ug42}

Here is our definition of Kuranishi cohomology~\cite[\S 4.4]{Joyc1}.

\begin{dfn} Let $Y$ be an orbifold without boundary. Consider
triples $(X,\bs f,\bs C)$, where $X$ is a compact Kuranishi space,
$\bs f:X\ra Y$ is a strong submersion with $(X,\bs f)$ {\it
cooriented}, as in \S\ref{ug27}, and $\bs C$ is {\it co-gauge-fixing
data\/} for $(X,\bs f)$, as in \S\ref{ug41}. Write $[X,\bs f,\bs C]$
for the isomorphism class of $(X,\bs f,\bs C)$ under isomorphisms
$(\bs a,\bs b):(X,\bs f,\bs C)\ra(\ti X,\bs{\ti f},\bs{\ti C})$,
where $\bs a$ must identify the coorientations of~$(X,\bs f),(\ti
X,\bs{\ti f})$, and $\bs b$ lifts $\bs a$ to the Kuranishi
neighbourhoods $(V^i,\ldots,\psi^i),(\ti V^i,\ldots,\ti\psi^i)$
in~$\bs C,\bs{\ti C}$.

Let $R$ be a $\Q$-algebra. For $k\in\Z$, define $KC^k(Y;R)$ to be
the $R$-module of finite $R$-linear combinations of isomorphism
classes $[X,\bs f,\bs C]$ for which $\vdim X=\dim Y-k$, with the
analogues of relations Definition \ref{ug3def}(i)--(iii), replacing
gauge-fixing data $\bs G$ by co-gauge-fixing data $\bs C$. Elements
of $KC^k(Y;R)$ are called {\it Kuranishi cochains}. Define
$\d:KC^k(Y;R)\ra KC^{k+1}(Y;R)$ by
\e
\d:\ts\sum_{a\in A}\rho_a[X_a,\bs f_a,\bs C_a]\longmapsto
\ts\sum_{a\in A}\rho_a[\pd X_a,\bs f_a\vert_{\pd X_a},\bs C_a
\vert_{\pd X_a}].
\label{ug4eq1}
\e
As in Definition \ref{ug3def} we have $\d\ci\d=0$. Define the {\it
Kuranishi cohomology groups} $KH^k(Y;R)$ of $Y$ for $k\in\Z$ to be
\begin{equation*}
KH^k(Y;R)=\frac{\Ker\bigl(\d:KC^k(Y;R)\ra KC^{k+1}(Y;R)\bigr)}{
\Im\bigl(\d:KC^{k-1}(Y;R)\ra KC^k(Y;R)\bigr)}\,.
\end{equation*}

Let $Y,Z$ be orbifolds without boundary, and $h:Y\ra Z$ be a smooth,
proper map. Define the {\it pullback\/} $h^*:KC^k(Z;R)\ra KC^k(Y;R)$
by
\e
h^*:\ts\sum_{a\in A}\rho_a\bigl[X_a,\bs f_a,\bs C_a\bigr]\longmapsto
\ts\sum_{a\in A}\rho_a [Y\t_{h,Z,\bs f_a}X_a,\bs\pi_Y,h^*(\bs C_a)],
\label{ug4eq2}
\e
where $h^*(\bs C_a)$ is as in Theorem \ref{ug4thm1}(f), and the
coorientation for $(X_a,\bs f_a)$ pulls back to a natural
coorientation for $(Y\t_ZX_a,\bs\pi_Y)$. These $h^*:KC^k(Z;R)\ra
KC^k(Y;R)$ satisfy $h^*\ci\d=\d\ci h^*$, so they induce morphisms of
cohomology groups $h^*:KH^k(Z;R)\ra KH^k(Y;R)$. Pullbacks are
functorial, that is, $(g\ci h)^*=h^*\ci g^*$, on both cochains and
cohomology.
\label{ug4def1}
\end{dfn}

Here for simplicity we restrict to orbifolds $Y$ {\it without
boundary}. When $\pd Y\ne\es$, the definition of $KH^*(Y;R)$ is more
complicated \cite[\S 4.5]{Joyc1}: if $\bs f:X\ra Y$ is a strong
submersion then we must split $\pd X=\pd_+^{\bs f}X\amalg\pd_-^{\bs
f}X$, where roughly speaking $\pd_+^{\bs f}X$ is the component of
$\pd X$ lying over $Y^\ci$, and $\pd_-^{\bs f}X$ the component of
$\pd X$ lying over $\pd Y$. Then $\bs f\vert_{\pd_+^{\bs f}X}$ is a
strong submersion $\bs f_+:\pd_+^{\bs f}X\ra Y$, and $\bs
f\vert_{\pd_-^{\bs f}X}=\io\ci\bs f_-$, where $\bs f_-:\pd_-^{\bs
f}X\ra\pd Y$ is a strong submersion, and $\io:\pd Y\ra Y$ is the
natural immersion. In \eq{ug4eq1} we must replace $[\pd X_a,\bs
f_a\vert_{\pd X_a},\bs C_a \vert_{\pd X_a}]$ by~$[\pd_+^{\bs
f_a}X_a,\bs f_{a,+},\bs C_a \vert_{\pd_+^{\bs f_a}X_a}]$.

In \cite[\S 4.7]{Joyc1} we define {\it cup\/} and {\it cap
products}.

\begin{dfn} Let $Y$ be an orbifold without boundary, and $R$ a
$\Q$-algebra. Define the {\it cup product\/} $\cup:KC^k(Y;R)\t
KC^l(Y;R)\ra KC^{k+l}(Y;R)$ by
\begin{equation*}
[X,\bs f,\bs C]\cup[\ti X,\bs{\ti f},\bs{\ti C}]=
\smash{\bigl[X\t_{\bs f,Y,\bs{\ti f}}\ti X,\bs\pi_Y,\bs C\t_Y
\bs{\ti C}\bigr]},
\end{equation*}
extended $R$-bilinearly. Here $\bs\pi_Y:X\t_{\smash{\bs f,Y,\bs{\ti
f}}}\ti X\ra Y$ is the projection from the fibre product, which is a
strong submersion as both $\bs f,\bs{\ti f}$ are, and $\bs
C\t_Y\bs{\ti C}$ is as in Theorem \ref{ug4thm1}(g). Then $\cup$
takes relations (i)--(iii) in both $KC^k(Y;R)$ and $KC^l(Y;R)$ to
the same relations in $KC^{k+l}(Y;R)$. Thus $\cup$ is well-defined.

Theorem \ref{ug4thm1}(g) and \eq{ug2eq13}--\eq{ug2eq15} imply that
for $\ga\in KC^k(Y;R)$, $\de\in KC^l(Y;R)$ and $\ep\in KC^m(Y;R)$ we
have
\begin{gather}
\ga\cup\de=(-1)^{kl}\de\cup\ga,
\label{ug4eq3}\\
\d(\ga\cup\de)=(\d\ga)\cup\de+(-1)^k\ga\cup(\d\de)\;\>\text{and}\;\>
(\ga\cup\de)\cup\ep=\ga\cup(\de\cup\ep).
\label{ug4eq4}
\end{gather}
Therefore in the usual way $\cup$ induces an associative,
supercommutative product $\cup:KH^k(Y;R)\t KH^l(Y;R)\ra
KH^{k+l}(Y;R)$ given for $\ga\in KC^k(Y;R)$ and $\de\in KC^l(Y;R)$
with $\d\ga=\d\de=0$ by
\begin{equation*}
(\ga+\Im\d_{k-1})\cup(\de+\Im\d_{l-1})=(\ga\cup\de)+\Im\d_{k+l-1}.
\end{equation*}

Now suppose $Y$ is compact. Then $\id_Y:Y\ra Y$ is a (strong)
submersion, with a trivial coorientation giving the positive sign to
the zero vector bundle over $Y$. One can define natural
co-gauge-fixing data $\bs C_Y$ for $(Y,\id_Y)$ such that
$[Y,\id_Y,\bs C_Y]\in KC^0(Y;R)$, with $\d[Y,\id_Y,\bs C_Y]=0$, and
for all $[X,\bs f,\bs C]\in KC^k(Y;R)$ we have
\begin{equation*}
[Y,\id_Y,\bs C_Y]\cup [X,\bs f,\bs C]=[X,\bs f,\bs
C]\cup[Y,\id_Y,\bs C_Y]=[X,\bs f,\bs C].
\end{equation*}
Thus $[Y,\id_Y,\bs C_Y]$ is the identity for $\cup$, at the cochain
level. Passing to cohomology, $\bigl[[Y,\id_Y,\bs C_Y]\bigr]\in
KH^0(Y;R)$ is the identity for $\cup$ in $KH^*(Y;R)$. We call
$\bigl[[Y,\id_Y,\bs C_Y]\bigr]$ the {\it fundamental class\/}
of~$Y$.

Define the {\it cap product\/} $\cap:KC_k(Y;R)\t KC^l(Y;R)\ra
KC_{k-l}(Y;R)$ by
\begin{equation*}
\smash{[X,\bs f,\bs G]\cap[\ti X,\bs{\ti f},\bs{\ti C}]=\bigl[X
\t_{\bs f,Y,\bs{\ti f}}\ti X,\bs\pi_Y,\bs G\t_Y\bs{\ti C}\bigr]},
\end{equation*}
extended $R$-bilinearly, with $\bs G\t_Y\bs{\ti C}$ as in Theorem
\ref{ug4thm1}(h). For $\ga\in KC_k(Y;R)$ and $\de,\ep\in KC^*(Y;R)$,
the analogue of \eq{ug4eq4}, using Theorem \ref{ug4thm1}(h) and
\eq{ug2eq16},~is
\begin{equation*}
\pd(\ga\cap\de)=(\pd\ga)\cap\de+(-1)^{\dim Y-k}\ga\cup(\d\de),
\quad(\ga\cap\de)\cap\ep=\ga\cap(\de\cup\ep).
\end{equation*}
If also $Y$ is compact then $\ga\cap[Y,\id_Y,\bs C_Y]=\ga$. Thus
$\cap$ induces a {\it cap product\/} $\cap:KH_k(Y;R)\t KH^l(Y;R)\ra
KH_{k-l}(Y;R)$. These products $\cap$ make Kuranishi chains and
homology into {\it modules} over Kuranishi cochains and cohomology.

Let $Y,Z$ be orbifolds without boundary, and $h:Y\ra Z$ a smooth,
proper map. Then \cite[Prop.~3.33]{Joyc1} implies that pullbacks
$h^*$ and pushforwards $h_*$ are compatible with $\cup,\cap$ on
(co)chains, in the sense that if $\al\in KC_*(Y;R)$ and $\be,\ga\in
KC^*(Z;R)$ then
\e
h^*(\be\cup\ga)=h^*(\be)\cup h^*(\ga) \quad\text{and}\quad
h_*(\al\cap h^*(\be))=h_*(\al)\cap\be.
\label{ug4eq5}
\e
Since $\cup,\cap,h^*,h_*$ are compatible with $\d,\pd$, passing to
(co)homology shows that \eq{ug4eq5} also holds for $\al\in
KH_*(Y;R)$ and $\be,\ga\in KH^*(Z;R)$. If $Z$ is compact then $Y$
is, with $h^*\bigl([Z,\id_Z,\bs C_Z]\bigr)=[Y,\id_Y,\bs C_Y]$ in
$KC^0(Y;R)$ and $h^*\bigl(\bigl[[Z,\id_Z,\bs C_Z]\bigr]\bigr)\ab=
\bigl[[Y,\id_Y,\bs C_Y]\bigr]$ in~$KH^0(Y;R)$.

To summarize: Kuranishi cochains $KC^*(Y;R)$ form a {\it
supercommutative, associative, differential graded\/ $R$-algebra},
and Kuranishi cohomology $KH^*(Y;R)$ is a {\it supercommutative,
associative, graded\/ $R$-algebra}. These algebras are {\it with
identity\/} if $Y$ is compact without boundary, and {\it without
identity\/} otherwise. Pullbacks $h^*$ induce {\it algebra
morphisms\/} on both cochains and cohomology. Kuranishi chains
$KC_*(Y;R)$ are a {\it graded module\/} over $KC^*(Y;R)$, and
Kuranishi homology $KH_*(Y;R)$ is a {\it graded module\/}
over~$KH^*(Y;R)$.
\label{ug4def2}
\end{dfn}

\subsection{Poincar\'e duality, and isomorphism with $H^*_\cs(Y;R)$}
\label{ug43}

Suppose $Y$ is an oriented manifold, of dimension $n$, without
boundary, and not necessarily compact, and $R$ is a commutative
ring. Then as in Bredon \cite[\S VI.9]{Bred} there are {\it
Poincar\'e duality isomorphisms}
\e
\Pd:H^k_\cs(Y;R)\longra H_{n-k}^\rsi(Y;R)
\label{ug4eq6}
\e
between compactly-supported cohomology, and singular homology. If
$Y$ is also {\it compact\/} then it has a {\it fundamental class}
$[Y]\in H_n(Y;R)$, and we can write the Poincar\'e duality map $\Pd$
of \eq{ug4eq6} in terms of the cap product by $\Pd(\al)=[Y]\cap\al$
for $\al\in H^k_\cs(Y;R)$. Satake \cite[Th.~3]{Sata} showed that
Poincar\'e duality isomorphisms \eq{ug4eq6} exist when $Y$ is an
oriented orbifold without boundary and $R$ is a $\Q$-{\it algebra}.

Let $Y$ be an orbifold of dimension $n$ without boundary, and $R$ a
$\Q$-algebra. We wish to construct an isomorphism $\Pi_\cs^\Kch:
H^*_\cs(Y;R)\ra KH^*(Y;R)$ from compactly-supported cohomology to
Kuranishi cohomology. In the case in which $Y$ is oriented we will
define $\Pi_\cs^\Kch$ to be the composition
\e
\smash{\xymatrix@C=25pt{ H^k_\cs(Y;R) \ar[r]^(0.45){\Pd} &
H_{n-k}^\rsi(Y;R) \ar[r]^(0.45){\Pi_\rsi^\Kh} & KH_{n-k}(Y;R)
\ar[r]^{\Pi^\Kch_\Kh} & KH^k(Y;R), }}
\label{ug4eq7}
\e
where the isomorphism $\Pd$ is as in \eq{ug4eq6}, and $\Pi_\rsi^\Kh$
is as in \eq{ug3eq2} and is an isomorphism by Theorem \ref{ug3thm2},
and $\Pi^\Kch_\Kh$ is an isomorphism between Kuranishi (co)homology,
with inverse~$\Pi_\Kch^\Kh:KH^k(Y;R)\ra KH_{n-k}(Y;R)$.

These $\Pi_\Kch^\Kh,\Pi^\Kch_\Kh$ are defined in
\cite[Def.~4.14]{Joyc1}. At the (co)chain level, we define
$\Pi_\Kch^\Kh:KC^k(Y;R)\ra KC_{n-k}(Y;R)$ by $\Pi_\Kch^\Kh:[X,\bs
f,\bs C]\mapsto [X,\bs f,\bs G_{\bs C}]$, where $\bs G_{\bs C}$ is
gauge-fixing data for $(X,\bs f)$ constructed from the
co-gauge-fixing data $\bs C$ is a functorial way, and as $\bs f:X\ra
Y$ is cooriented and $Y$ is oriented, we obtain an orientation for
$X$ as in \S\ref{ug27}. Then $\pd\ci
\Pi_\Kch^\Kh=\Pi_\Kch^\Kh\ci\d$, so they induce morphisms
$\Pi_\Kch^\Kh:KH^k(Y;R)\ra KH_{n-k}(Y;R)$ in (co)homology.

For $\Pi^\Kch_\Kh$ the story is more complicated. To define
$\Pi^\Kch_\Kh:KC_{n-k}(Y;R)\ra KC^k(Y;R)$ we cannot simply map
$[X,\bs f,\bs G]\mapsto[X,\bs f,\bs C_{\bs G}]$ for some
co-gauge-fixing data $\bs C_{\bs G}$ constructed from $\bs G$, since
$\bs f$ need only be strongly smooth for $[X,\bs f,\bs G]\in
KC_{n-k}(Y;R)$, but $\bs f$ must be a strong submersion for $[X,\bs
f,\bs C_{\bs G}]\in KC^k(Y;R)$. Instead, we define $\Pi^\Kch_\Kh:
KC_{n-k}(Y;R)\ra KC^k(Y;R)$ by $\Pi^\Kch_\Kh:[X,\bs f,\bs G]\mapsto
[X^Y,\bs f{}^Y,\bs C_{\bs G}^Y]$. Here $X^Y$ is $X$ equipped with an
{\it alternative Kuranishi structure}, which roughly speaking adds
copies of $\bs f^*(TY)$ to both the tangent bundle and obstruction
bundle of $X$. Also $\bs f{}^Y:X^Y\ra Y$ is a lift of $\bs f$ to
$X^Y$, which is a strong submersion, and $\bs C_{\bs G}^Y$ is
co-gauge-fixing data for $(X^Y,\bs f{}^Y)$ constructed in a
functorial way from $\bs G$. Then $\d\ci\Pi^\Kch_\Kh=
\Pi^\Kch_\Kh\ci\pd$, so they induce morphisms
$\Pi^\Kch_\Kh:KH_{n-k}(Y;R)\ra KH^k(Y;R)$ on (co)homology.

In \cite[Th.~4.15]{Joyc1} we show that these morphisms
$\Pi_\Kch^\Kh,\Pi^\Kch_\Kh$ on Kuranishi (co)homology are inverse.
Thus the third morphism $\Pi^\Kch_\Kh$ in \eq{ug4eq7} is an
isomorphism, so the composition $\Pi_\cs^\Kch:H^k_\cs(Y;R)\ra
KH^k(Y;R)$ is an isomorphism. Changing the orientation of $Y$
changes the sign of $\Pd,\Pi^\Kch_\Kh$, and so does not change
$\Pi_\cs^\Kch$. If $Y$ is not orientable we can make a similar
argument using homology groups $H_{n-k}^\rsi(Y;O\t_{\{\pm 1\}}R)$,
$KH_{n-k}(Y;O\t_{\{\pm 1\}}R)$ twisted by the principal
$\Z_2$-bundle $O$ of orientations on $Y$. Thus we
prove~\cite[Cor.~4.17]{Joyc1}:

\begin{thm} Let\/ $Y$ be an orbifold without boundary, and\/ $R$ a
$\Q$-algebra. Then there are natural isomorphisms\/ $\Pi_\cs^\Kch:
H^k_\cs(Y;R)\ra KH^k(Y;R)$ for $k\ge 0,$ and\/ $KH^k(Y;R)=0$
when~$k<0$.
\label{ug4thm2}
\end{thm}

In \cite[\S 4.5]{Joyc1} the theorem is extended to $Y$ with
boundary, going via relative homology $H_*^\rsi(Y,\pd Y;R),
KH_*(Y,\pd Y;R)$. In \cite[Th.~4.34]{Joyc1} we show that the
isomorphisms $\Pi_\cs^\Kch:H^*_\cs(Y;R)\ra KH^*(Y;R)$ and
$\Pi_\rsi^\Kh:H_*^\rsi(Y;R)\ra KH_*(Y;R)$ in Theorems \ref{ug3thm2}
and \ref{ug4thm2} identify the cup and cap products $\cup,\cap$ on
$H^*_\cs(Y;R),H_*^\rsi(Y;R)$ with those on~$KH^*(Y;R),KH_*(Y;R)$.

\section{Kuranishi bordism and cobordism}
\label{ug5}

We now summarize parts of \cite[\S 5]{Joyc1} on Kuranishi
(co)bordism. They are based on the classical bordism theory
introduced by Atiyah \cite{Atiy}. In fact \cite[\S 5]{Joyc1} studies
five different kinds of Kuranishi (co)bordism, but we discuss only
one.

\subsection{Classical bordism and cobordism groups}
\label{ug51}

Bordism groups were introduced by Atiyah \cite{Atiy}, and Connor
\cite[\S I]{Conn} gives a good introduction. Our definition is not
standard, but fits in with~\S\ref{ug52}.

\begin{dfn} Let $Y$ be an orbifold without boundary. Consider pairs
$(X,f)$, where $X$ is a compact, oriented manifold without boundary
or corners, not necessarily connected, and $f:X\ra Y$ is a smooth
map. An {\it isomorphism\/} between two such pairs $(X,f),(\ti X,\ti
f)$ is an orientation-preserving diffeomorphism $i:X\ra\ti X$ with
$f=\ti f\ci i$. Write $[X,f]$ for the isomorphism class of~$(X,f)$.

Let $R$ be a commutative ring. For each $k\ge 0$, define the $k^{\it
th}$ {\it bordism group\/} $B_k(Y;R)$ of $Y$ with coefficients in
$R$ to be the $R$-module of finite $R$-linear combinations of
isomorphism classes $[X,f]$ for which $\dim X=k$, with the
relations:
\begin{itemize}
\setlength{\itemsep}{0pt}
\setlength{\parsep}{0pt}
\item[(i)] $[X,f]+[X',f']=[X\amalg X',f\amalg f']$ for all
classes $[X,f],[X',f']$; and
\item[(ii)] Suppose $Z$ is a compact, oriented $(k\!+\!1)$-manifold
with boundary but without (g-)corners, and $g:Z\ra Y$ is smooth.
Then~$[\pd Z,g\vert_{\pd Z}]=0$.
\end{itemize}
\label{ug5def1}
\end{dfn}

Here is how this definition relates to those in \cite{Atiy,Conn}.
When $Y$ is a manifold and $R=\Z$, our $B_k(Y;\Z)$ is equivalent to
Connor's {\it differential bordism group\/} $D_k(Y)$, \cite[\S
I.9]{Conn}. Atiyah \cite[\S 2]{Atiy} and Connor \cite[\S I.4]{Conn}
also define {\it bordism groups\/} $MSO_k(Y)$ as for $B_k(Y;\Z)$
above, but only requiring $f:X\ra Y$ to be continuous, not smooth.
Connor \cite[Th.~I.9.1]{Conn} shows that when $Y$ is a manifold, the
natural projection $D_k(Y)\ra MSO_k(Y)$ is an isomorphism.

As in \cite[\S I.5]{Conn}, bordism is a {\it generalized homology
theory}, that is, it satisfies all the Eilenberg--Steenrod axioms
for a homology theory except the dimension axiom. The bordism groups
of a point $MSO_*({\rm pt})$ are known, \cite[\S I.2]{Conn}. This
gives some information on bordism groups of general spaces $Y$: for
any generalized homology theory $GH_*(Y)$, there is a spectral
sequence from the singular homology $H^\rsi_* \bigl(Y;GH_*({\rm
pt})\bigr)$ of $Y$ with coefficients in $GH_*({\rm pt})$ converging
to $GH_*(Y)$, so that $GH_*({\cal S}^n)\cong H^\rsi_*\bigl({\cal
S}^n;GH_*({\rm pt})\bigr)$, for instance.

Atiyah \cite{Atiy} and Connor \cite[\S 13]{Conn} also define {\it
cobordism groups} $MSO^k(Y)$ for $k\in\Z$, which are a {\it
generalized cohomology theory\/} dual to bordism $MSO_k(Y)$. There
is a natural product $\cup$ on $MSO^*(Y)$, making it into a
supercommutative ring. If $Y$ is a compact, oriented $n$-manifold
without boundary then \cite[Th.~3.6]{Atiy}, \cite[Th.~13.4]{Conn}
there are canonical Poincar\'e duality isomorphisms
\e
MSO^k(Y)\cong MSO_{n-k}(Y)\quad\text{for $k\in\Z$.}
\label{ug5eq1}
\e
The definition of $MSO^*(Y)$ uses homotopy theory, direct limits of
$k$-fold suspensions, and classifying spaces. There does not seem to
be a satisfactory differential-geometric definition of cobordism
groups parallel to Definition~\ref{ug5def1}.

\subsection{Kuranishi bordism and cobordism groups}
\label{ug52}

Motivated by \S\ref{ug51}, following \cite[\S 5.2]{Joyc1} we define:

\begin{dfn} Let $Y$ be an orbifold. Consider pairs $(X,\bs f)$,
where $X$ is a compact, oriented Kuranishi space without boundary,
and $\bs f:X\ra Y$ is strongly smooth. An {\it isomorphism\/}
between two pairs $(X,\bs f),(\ti X,\bs{\ti f})$ is an
orientation-preserving strong diffeomorphism $\bs i:X\ra\ti X$ with
$\bs f=\bs{\ti f}\ci\bs i$. Write $[X,\bs f]$ for the isomorphism
class of~$(X,\bs f)$.

Let $R$ be a commutative ring. For each $k\in\Z$, define the $k^{\it
th}$ {\it Kuranishi bordism group\/} $KB_k(Y;R)$ of $Y$ with
coefficients in $R$ to be the $R$-module of finite $R$-linear
combinations of isomorphism classes $[X,\bs f]$ for which $\vdim
X=k$, with the relations:
\begin{itemize}
\setlength{\itemsep}{0pt}
\setlength{\parsep}{0pt}
\item[(i)] $[X,\bs f]+[X',\bs f']=[X\amalg X',\bs f\amalg
\bs f']$ for all classes $[X,\bs f],[X',\bs f']$; and
\item[(ii)] Suppose $W$ is a compact, oriented Kuranishi space with
boundary but without (g-)corners, with $\vdim W=k+1$, and $\bs
e:W\ra Y$ is strongly smooth. Then~$[\pd W,\bs e\vert_{\pd W}]=0$.
\end{itemize}
Elements of $KB_k(Y;R)$ will be called {\it Kuranishi bordism
classes}.

Let $h:Y\ra Z$ be a smooth map of orbifolds. Define the {\it
pushforward\/} $h_*:KB_k(Y;R)\ra KB_k(Z;R)$ by $h_*:\sum_{a\in
A}\rho_a[X_a,\bs f_a]\mapsto\sum_{a\in A}\rho_a[X_a,h\ci\bs f_a]$.
This takes relations (i),(ii) in $KB_k(Y;R)$ to (i),(ii) in
$KB_k(Z;R)$, and so is well-defined. Pushforward is functorial, that
is,~$(g\ci h)_*=g_*\ci h_*$.
\label{ug5def2}
\end{dfn}

Now Kuranishi bordism $KB_*(Y;R)$ is like Kuranishi homology
$KH_*(Y;R)$ in \S\ref{ug32}, but using Kuranishi spaces $X$ without
boundary, and omitting gauge-fixing data $\bs G$. Thus it seems
natural to define Kuranishi cobordism $KB^*(Y;R)$ by modifying the
definition of Kuranishi cohomology $KH^*(Y;R)$ in \S\ref{ug42} in
the same way, following~\cite[\S 5.4--\S 5.5]{Joyc1}.

\begin{dfn} Let $Y$ be an orbifold without boundary. Consider pairs
$(X,\bs f)$, where $X$ is a compact Kuranishi space without
boundary, and $\bs f:X\ra Y$ is a cooriented strong submersion. An
{\it isomorphism\/} between two pairs $(X,\bs f),(\ti X,\bs{\ti f})$
is a coorientation-preserving strong diffeomorphism $\bs i:X\ra\ti
X$ with $\bs f=\bs{\ti f}\ci\bs i$. Write $[X,\bs f]$ for the
isomorphism class of~$(X,\bs f)$.

Let $R$ be a commutative ring. For each $k\in\Z$, define the $k^{\it
th}$ {\it Kuranishi cobordism group\/} $KB^k(Y;R)$ of $Y$ with
coefficients in $R$ to be the $R$-module of finite $R$-linear
combinations of isomorphism classes $[X,\bs f]$ for which $\vdim
X=\dim Y-k$, with the relations:
\begin{itemize}
\setlength{\itemsep}{0pt}
\setlength{\parsep}{0pt}
\item[(i)] $[X,\bs f]+[X',\bs f']=[X\amalg X',\bs f\amalg
\bs f']$ for all classes $[X,\bs f],[X',\bs f']$; and
\item[(ii)] Suppose $W$ is a compact Kuranishi space with
boundary but without \hbox{(g-)}\ab corners, with $\vdim W=\dim
Y-k+1$, and $\bs e:W\ra Y$ is a cooriented strong submersion. Then
$\bs e\vert_{\pd W}:\pd W\ra Y$ is a cooriented strong submersion,
and we impose the relation $[\pd W,\bs e\vert_{\pd W}]=0$
in~$KB^k(Y;R)$.
\end{itemize}
Elements of $KB^k(Y;R)$ will be called {\it Kuranishi cobordism
classes}.

Define the {\it cup product\/} $\cup:KB^k(Y;R)\t KB^l(Y;R)\ra
KB^{k+l}(Y;R)$ by
\e
\raisebox{-4pt}{\begin{Large}$\displaystyle\Bigl[$\end{Large}}
\sum_{a\in A}\rho_a\bigl[X_a,\bs f_a\bigr]
\raisebox{-4pt}{\begin{Large}$\displaystyle\Bigr]$\end{Large}}\!\cup\!
\raisebox{-4pt}{\begin{Large}$\displaystyle\Bigl[$\end{Large}}
\sum_{b\in B}\si_b\bigl[\ti X_b,\bs{\ti f}_b\bigr]
\raisebox{-4pt}{\begin{Large}$\displaystyle\Bigr]$\end{Large}}\!=\!\!
\sum_{a\in A,\; b\in B\!\!\!\!\!\!}\!\!\rho_a\si_b
\bigl[X_a\t_{\bs f_a,Y,\bs{\ti f}_b}\ti X_b,\bs\pi_Y\bigr],
\label{ug5eq2}
\e
for $A,B$ finite and $\rho_a,\si_b\in R$. The coorientations on
$(X_a,\bs f_a)$ and $(\ti X_b,\bs{\ti f}_b)$ induce a coorientation
on $(X_a\t_Y\ti X_b,\bs\pi_Y)$ as in \S\ref{ug27}. Similarly, define
the {\it cap product\/} $\cap:KB_k(Y;R)\t KB^l(Y;R)\ra
KB_{k-l}(Y;R)$ by the same formula \eq{ug5eq2}, where now $[X_a,\bs
f_a]\in KB_k(Y;R)$ so that $\bs f_a$ is strongly smooth and $X_a$
oriented, and the orientation on $X_a$ and coorientation for
$\bs{\ti f}_b$ combine to give an orientation for~$X_a\t_Y\ti X_b$.

One can show that $\cup,\cap$ are well-defined, that $\cup$ is
associative and supercommutative, and that
$(\ga\cap\de)\cap\ep=\ga\cap(\de\cup\ep)$ for $\ga\in KB_*(Y;R)$ and
$\de,\ep\in KB^*(Y;R)$. If $Y$ is also {\it compact\/} then using
the trivial coorientation for $\id_Y:Y\ra Y$, we have $[Y,\id_Y]\in
KB^0(Y;R)$, which is the {\it identity\/} for $\cup$ and $\cap$.
Thus, $KB^*(Y;R)$ is a {\it graded, supercommutative, associative
$R$-algebra, with identity\/} if $Y$ is compact, and {\it without
identity\/} otherwise, and $\cap$ makes $KB_*(Y;R)$ into a module
over~$KB^*(Y;R)$.

Let $Y,Z$ be orbifolds without boundary, and $h:Y\ra Z$ a smooth,
proper map. Motivated by \eq{ug4eq2}, define the {\it pullback\/}
$h^*:KB^k(Z;R)\ra KB^k(Y;R)$ by $h^*:\sum_{a\in A}\rho_a[X_a,\bs
f_a]\mapsto\sum_{a\in A}\rho_a [Y\t_{h,Z,\bs f_a}X_a,\bs\pi_Y]$.
This takes relations (i),(ii) in $KB^k(Z;R)$ to (i),(ii) in
$KB^k(Y;R)$, and so is well-defined. Pullbacks are functorial,
$(g\ci h)^*=h^*\ci g^*$. The cup and cap products are compatible
with pullbacks and pushforwards, as in~\eq{ug4eq5}.
\label{ug5def3}
\end{dfn}

\subsection{Morphisms to and from Kuranishi (co)bordism}
\label{ug53}

In \cite[\S 5.3--\S 5.4]{Joyc1} we define morphisms between these
groups.

\begin{dfn} Let $Y$ be an orbifold, and $R$ a commutative ring.
Define morphisms $\Pi_\bo^\Kb:B_k(Y;R)\ra KB_k(Y;R)$ for $k\ge 0$ by
$\Pi_\bo^\Kb:\sum_{a\in A}\rho_a[X_a,f_a]\mapsto
\sum_{a\in A}\rho_a[X_a,f_a]$, interpreting the manifold $X_a$ as a
Kuranishi space, and the smooth map $f_a:X_a\ra Y$ as strongly
smooth.

Define morphisms $\Pi_\Kb^\Kh:KB_k(Y;R)\ra KH_k(Y;R\ot_\Z\Q)$ for
$k\in\Z$ by $\Pi_\Kb^\Kh:\sum_{a\in A}\rho_a\bigl[X_a,\bs
f_a\bigr]\mapsto\bigl[\ts\sum_{a\in A}\pi(\rho_a)[X_a,\ab\bs f_a,\bs
G_a]\bigr]$, where $\bs G_a$ is some choice of gauge-fixing data for
$(X_a,\bs f_a)$, which exists by Theorem \ref{ug3thm1}(a), and
$\pi:R\ra R\ot_\Z\Q$ is the natural morphism. Using Theorem
\ref{ug3thm1}(e) over $[0,1]\t X_a$ one can show that $\Pi_\Kb^\Kh$
is independent of the choice of $\bs G_a$, and is well-defined.

Similarly, define morphisms $\Pi_\Kcb^\Kch:KB^k(Y;R)\ra
KH^k(Y;R\ot_\Z\Q)$ for $k\in\Z$ by $\Pi_\Kb^\Kh:\sum_{a\in
A}\rho_a\bigl[X_a,\bs f_a\bigr]\mapsto\bigl[\ts\sum_{a\in
A}\pi(\rho_a)[X_a,\ab\bs f_a,\bs C_a]\bigr]$, where $\bs C_a$ is
some choice of co-gauge-fixing data for $(X_a,\bs f_a)$.

These $\Pi_\Kcb^\Kch,\Pi_\Kb^\Kh$ take cup and cap products
$\cup,\cap$ on $KB^*,KB_*(Y;R)$ to $\cup,\cap$ on $KH^*,KH_*
(Y;R\ot_\Z\Q)$, and if $Y$ is compact, they take the identity
$[Y,\id_Y]\in KB^0(Y;R)$ to the identity $\bigl[[Y,\id_Y,\bs
C_Y]\bigr]\in KH^0(Y;R\ot_\Z\Q)$.
\label{ug5def4}
\end{dfn}

Consider the sequence of morphisms
\begin{equation*}
\smash{\xymatrix@C=25pt{ B_k(Y;R) \ar[r]^(0.45){\Pi_\bo^\Kb} &
KB_k(Y;R) \ar[r]^(0.4){\Pi_\Kb^\Kh} & KH_k(Y;R\ot_\Z\Q)
\ar[r]^(0.55){(\Pi_\rsi^\Kh)^{\smash{-1}}} &
H^\rsi_k(Y;R\ot_\Z\Q),}}
\end{equation*}
where $(\Pi_\rsi^\Kh)^{-1}$ exists by Theorem \ref{ug3thm2}. The
composition is the natural map $B_k(Y;R)\ra H^\rsi_k(Y;R\ot_\Z\Q)$
taking $[X,f]\mapsto f_*([X])$. Thus we find:

\begin{cor} Let\/ $Y$ be an orbifold, and\/ $R$ a commutative ring.
Then $KB_k(Y;\ab R)$ is at least as large as the image of\/
$B_k(Y;R)$ in\/~$H_k^\rsi(Y;R\ot_\Z\Q)$.
\label{ug5cor}
\end{cor}

We will see in \S\ref{ug54} that $KB_*(Y;R)$ is actually very large.

The Poincar\'e duality story for Kuranishi (co)homology in
\S\ref{ug43} has an analogue for Kuranishi (co)bordism, as in
\cite[\S 5.4]{Joyc1}. Let $Y$ be an oriented $n$-orbifold without
boundary, and $R$ a commutative ring. Define $R$-module morphisms
$\Pi_\Kcb^\Kb:KB^k(Y;R)\ra KB_{n-k}(Y;R)$ for $k\in\Z$ by
$\Pi_\Kcb^\Kb:\sum_{a\in A}\rho_a[X_a,\bs f_a]\mapsto\sum_{a\in
A}\rho_a[X_a,\bs f_a]$, using the coorientation for $\bs f_a$ from
$[X_a,\bs f_a]\in KB^k(Y;R)$ and the orientation on $Y$ to determine
the orientation on $X_a$ for $[X_a,\bs f_a]\in KB_{n-k}(Y;R)$.

Define $\Pi^\Kcb_\Kb:KB_{n-k}(Y;R)\ra KB^k(Y;R)$ by
$\Pi^\Kcb_\Kb:\sum_{a\in A}\rho_a[X_a,\ab\bs f_a]\mapsto
\sum_{a\in A}\rho_a[X_a^Y,\bs f{}_a^Y]$, where $X_a^Y$ is $X_a$ with
an {\it alternative Kuranishi structure\/} as in \S\ref{ug43}, and
$\bs f_a{}^Y:X_a^Y\ra Y$ is a lift of $\bs f_a$ to $X_a^Y$, which is
a strong submersion. Then \cite[Th.~5.11]{Joyc1} shows that
$\Pi^\Kcb_\Kb$ and $\Pi_\Kcb^\Kb$ are inverses, so they are both
isomorphisms. Using this and ideas in \S\ref{ug51} including
\eq{ug5eq1}, if $Y$ is a compact manifold we can define a natural
morphism $\Pi_{\rm cb}^\Kcb:MSO^*(Y)\ra KB^*(Y;\Z)$, so Kuranishi
cobordism is a generalization of classical cobordism.

\subsection{How large are Kuranishi (co)bordism groups?}
\label{ug54}

Theorems \ref{ug3thm2} and \ref{ug4thm2} showed that Kuranishi
(co)homology are isomorphic to classical (compactly-supported)
(co)homology, so they are not new topological invariants. In
contrast, Kuranishi (co)bordism are not isomorphic to classical
(co)bordism, they are genuinely new topological invariants, so it is
interesting to ask what we can say about them. We now summarize the
ideas of \cite[\S 5.6--\S 5.7]{Joyc1}, which show that $KB_*(Y;R)$
and $KB^*(Y;R)$ are {\it very large\/} for any orbifold $Y$ and
commutative ring $R$ with $Y\ne\es$ and~$R\ot_\Z\Q\ne 0$.

One reason for this is that in a class $\sum_{a\in A}\rho_a[X_a,\bs
f_a]$ in $KB_k(Y;R)$ there is a lot of information stored in the
{\it orbifold strata\/} of $X_a$ for $a\in A$. We define these for
orbifolds,~\cite[Def.~5.15]{Joyc1}.

\begin{dfn} Let $\Ga$ be a finite group, and consider
(finite-dimensional) real representations $(W,\om)$ of $\Ga$, that
is, $W$ is a finite-dimensional real vector space and
$\om:\Ga\ra\Aut(W)$ is a group morphism. Call $(W,\om)$ a {\it
trivial representation} if $\om\equiv\id_W$, and a {\it nontrivial
representation} if $\Fix(\om(\Ga))=\{0\}$. Then every
$\Ga$-representation $(W,\om)$ has a unique decomposition $W=W^{\rm
t}\op W^{\rm nt}$ as the direct sum of a trivial representation
$(W^{\rm t},\om^{\rm t})$ and a nontrivial representation $(W^{\rm
nt},\om^{\rm nt})$, where~$W^{\rm t}=\Fix(\om(\Ga))$.

Now let $X$ be an $n$-orbifold, $\Ga$ be a finite group, and $\rho$
be an {\it isomorphism class of nontrivial\/ $\Ga$-representations}.
Each $p\in X$ has a {\it stabilizer group\/} $\Stab_X(p)$. The
tangent space $T_pX$ is an $n$-dimensional vector space with a
representation $\tau_p$ of $\Stab_X(p)$. Let $\la:\Ga\ra\Stab_X(p)$
be an injective group morphism, so that $\la(\Ga)$ is a subgroup of
$\Stab_X(p)$ isomorphic to $\Ga$. Hence
$\tau_p\ci\la:\Ga\ra\Aut(T_pX)$ is a $\Ga$-representation, and we
can split $T_pX=(T_pX)^{\rm t}\op(T_pX)^{\rm nt}$ into trivial and
nontrivial $\Ga$-representations, and form the isomorphism class
$\bigl[(T_pX)^{\rm nt},(\tau_p\ci\la)^{\rm nt}\bigr]$. As a set,
define the {\it orbifold stratum} $X^{\Ga,\rho}$ to be
\begin{align*}
X^{\Ga,\rho}=\bigl\{\Stab_X(p)\cdot(p,\la):\text{$p\in X$,
$\la:\Ga\ra\Stab_X(p)$ is an injective}& \\
\text{group morphism, $\bigl[(T_pX)^{\rm nt},(\tau_p\ci\la)^{\rm
nt}\bigr]=\rho$}&\bigr\},
\end{align*}
where $\Stab_X(p)$ acts on pairs $(p,\la)$ by $\si:(p,\la)\mapsto
(p,\la^\si)$, where $\la^\si:\Ga\ra\Stab_X(p)$ is given by
$\la^\si(\ga)=\si\la(\ga)\si^{-1}$. Define a map $\io^{\Ga,\rho}:
X^{\Ga,\rho}\ra X$ by~$\io^{\Ga,\rho}:\Stab_X(p)\cdot(p,\la)\mapsto
p$. Then \cite[Prop.~5.16]{Joyc1} shows that $X^{\Ga,\rho}$ is an
orbifold of dimension $n-\dim\rho$, and $\io^{\Ga,\rho}$ lifts to a
proper, finite immersion.
\label{ug5def5}
\end{dfn}

If $X$ is a Kuranishi space, there is a parallel definition
\cite[Def.~5.18]{Joyc1} of orbifold strata $X^{\Ga,\rho}$, which we
will not give. The most important difference is that $\rho$ is now a
{\it virtual nontrivial representation\/} of $\Ga$, that is, a
formal difference of nontrivial representations, so that
$\dim\rho\in\Z$ rather than $\dim\rho\in\N$. We find
\cite[Prop.~5.19]{Joyc1} that $X^{\Ga,\rho}$ is a Kuranishi space
with $\vdim X^{\Ga,\rho}=\vdim X-\dim\rho$, equipped with a proper,
finite, strongly smooth map~$\bs\io^{\Ga,\rho}:X^{\Ga,\rho}\ra X$.

We would like to define projections $\Pi^{\Ga,\rho}:KB_k(Y;R)\ra
KB_{k-\dim\rho}(Y;R)$ mapping $\Pi^{\Ga,\rho}:[X_a,\bs f_a]\ra
[X_a^{\Ga,\rho},\bs f_a\vert_{X_a^{\Ga,\rho}}]$. But there is a
problem: we need to define an {\it orientation\/} on
$X_a^{\Ga,\rho}$ from the orientation on $X_a$, and for general
$\Ga,\rho$ this may not be possible. To overcome this we suppose
$\md{\Ga}$ is odd, which implies that $\dim\rho$ is even for all
$\rho$, and there is then a consistent way to define orientations on
$X_a^{\Ga,\rho}$, and $\Pi^{\Ga,\rho}$ is well-defined.

Let $Y$ be a nonempty, connected orbifold. In \cite[\S 5.7]{Joyc1},
for each finite group $\Ga$ with $\md{\Ga}$ odd and all isomorphism
classes $\rho$ of virtual nontrivial representations of $\Ga$, we
construct a class $C^{\Ga,\rho}\in KB_{\dim\rho}(Y;\Z)$, such that
$\Pi_\Kb^\Kh\ci\Pi^{\Ga,\rho}(C^{\Ga,\rho})$ is nonzero in
$KH_0(Y;\Q)\cong H^\rsi_0(Y;\Q)\cong\Q$, and
$\Pi^{\De,\si}(C^{\Ga,\rho})=0$ if either $\md{\De}\ge\md{\Ga}$ and
$\De\not\cong\Ga$, or if $\De=\Ga$ and $\rho\ne\si$. It follows that
taken over all isomorphism classes of pairs $\Ga,\rho$, the classes
$C^{\Ga,\rho}\in KB_*(Y;\Z)$ are {\it linearly independent over\/}
$\Z$. Extending to an arbitrary commutative ring $R$, and using the
Poincar\'e duality ideas of \S\ref{ug53}, we deduce:

\begin{thm} Let\/ $Y$ be a nonempty orbifold, and\/ $R$ a
commutative ring with\/ $R\ot_\Z\Q\ne 0$. Then $KB_{2k}(Y;R)$ is
infinitely generated over $R$ for all\/ $k\in\Z$. If also $Y$ is
oriented of dimension $n$ then $KB^{n-2k}(Y;R)$ is infinitely
generated over $R$ for all\/~$k\in\Z$.
\label{ug5thm}
\end{thm}

Theorem \ref{ug5thm} supports the idea that Kuranishi bordism, (or
better, {\it almost complex Kuranishi bordism}, as in
\cite[Ch.~5]{Joyc1}) may be a useful tool for studying (closed)
Gromov--Witten invariants. In \cite[\S 6.2]{Joyc1} we define new
Gromov--Witten type invariants $[\oM_{g,m}(M,J,\be),\prod_i{\bf
ev}_i]$ in Kuranishi bordism $KB_*(M^m;\Z)$. Theorem \ref{ug5thm}
indicates that $KB_*(M^m;\Z)$ is very large, so that these new
invariants {\it contain a lot of information}, and that much of this
information has to do with the {\it orbifold strata\/} of the moduli
spaces $\oM_{g,m}(M,J,\be)$.

Also, these new invariants are defined in groups $KB_*(M^m;\Z)$ {\it
over $\Z$, not\/} $\Q$. When we project to Kuranishi homology or
singular homology to get conventional Gromov--Witten invariants, we
must work in homology over $\Q$. The reason we cannot work over $\Z$
is because of rational contributions from the orbifold strata of
$\oM_{g,m}(M,J,\be)$. Kuranishi bordism looks like a good framework
for describing these contributions, and so for understanding the
integrality properties of Gromov--Witten invariants, such as the
Gopakumar--Vafa Integrality Conjecture for Gromov--Witten invariants
of Calabi--Yau 3-folds. This is discussed in \cite[\S 6.3]{Joyc1},
and the author hopes to take it further in~\cite{Joyc4}.

\medskip

\noindent{\small\sc The Mathematical Institute, 24-29 St. Giles,
Oxford, OX1 3LB, U.K.}

\noindent{\small\sc E-mail: \tt joyce@maths.ox.ac.uk}


\begin{thebibliography}{99}

\bibitem{ALR} A. Adem, J. Leida and Y. Ruan, {\it Orbifolds and
Stringy Topology}, Cambridge Tracts in Math. 171, Cambridge
University Press, Cambridge, 2007.

\bibitem{AkJo} M. Akaho and D. Joyce, {\it Lagrangian Floer theory
using Kuranishi cohomology}, in preparation, 2009.

\bibitem{Atiy} M.F. Atiyah, {\it Bordism and cobordism}, Proc. Camb.
Phil. Soc. 57 (1961), 200--208.

\bibitem{Bred} G.E. Bredon, {\it Topology and Geometry}, Graduate
Texts in Math. 139, Springer-Verlag, New York, 1993.

\bibitem{ChSu} M. Chas and D. Sullivan, {\it String Topology},
math.DG/9911159, 1999. To appear in Annals of Mathematics.

\bibitem{Conn} P.E. Connor, {\it Differentiable Periodic Maps},
second edition, Springer Lecture Notes in Mathematics 738,
Springer-Verlag, Berlin, 1979.

\bibitem{Cost} K.J. Costello, {\it Topological conformal field
theories and Calabi-Yau categories}, Advances in Mathematics 210
(2007), 165--214. math.AG/0412149.

\bibitem{EES} T. Ekholm, J. Etynre and M. Sullivan, {\it The contact
homology of Legendrian submanifolds in $\R^{2n+1}$}, J. Diff. Geom.
71 (2005), 177--305.

\bibitem{EGH} Y. Eliashberg, A. Givental and H. Hofer, {\it
Introduction to Symplectic Field Theory}, Geom. Funct. Anal. 2000,
Special Volume, Part II, 560--673. math.SG/0010059.

\bibitem{FOOO} K. Fukaya, Y.-G. Oh, H. Ohta and K. Ono, {\it
Lagrangian intersection Floer theory -- anomaly and obstruction},
preprint, final(?) version, 2008. 1385 pages.

\bibitem{FuOn} K. Fukaya and K. Ono, {\it Arnold Conjecture and
Gromov--Witten invariant}, Topology 38 (1999), 933--1048.

\bibitem{Joyc1} D. Joyce, {\it Kuranishi homology and Kuranishi
cohomology}, arXiv:0707.3572, version 5, 2008. 290 pages.

\bibitem{Joyc2} D. Joyce, {\it A theory of open Gromov--Witten
invariants}, in preparation, 2009.

\bibitem{Joyc3} D. Joyce, {\it Kuranishi (co)homology and String
Topology}, in preparation, 2009.

\bibitem{Joyc4} D. Joyce, {\it Towards a symplectic proof of the
integrality conjecture for Gopakumar--Vafa invariants}, in
preparation, 2009.

\bibitem{LiTi} J. Li and G. Tian, {\it Virtual moduli cycles and
Gromov--Witten invariants of general symplectic manifolds}, pages
47--83 in R.J. Stern, editor, {\it Topics in symplectic
$4$-manifolds (Irvine, CA, 1996)}, International Press, Cambridge,
MA, 1998. alg-geom/9608032.

\bibitem{Ruan} Y. Ruan, {\it Virtual neighbourhoods and
pseudo-holomorphic curves}, Turkish J. Math. 23 (1999), 161--231.
alg-geom/96110121.

\bibitem{Sata} I. Satake, {\it On a generalization of the notion of
manifold}, Proc. Nat. Acad. Sci. U.S.A. 42 (1956), 359--363.

\bibitem{Sieb} B. Siebert, {\it Gromov--Witten invariants of
general symplectic manifolds}, dg-ga/9608005, 1996.

\end{thebibliography}
\end{document}